\documentstyle[a4,amstex,amscd,euscript,epsbox,amssymb]{amsart}
\numberwithin{equation}{section}
\newtheorem{thm}{Theorem}[section]
\newtheorem{cor}[thm]{Corollary}
\newtheorem{prop}[thm]{Proposition}
\newtheorem{lem}[thm]{Lemma}
\theoremstyle{remark}

\begin{document}
\newcommand{\dcp}{\,{\triangleright \kern-0.145em \triangleleft}}
\newcommand{\dcpP}{\,{\triangleright \kern-0.145em \triangleleft}_{\cal P}}
\newcommand{\tensE}{{\otimes}_{\frak E}}
\newcommand{\tenc}{{\intercal}}
\newcommand{\cu}{\varepsilon}
\newcommand{\emaru}{\stackrel{\scriptscriptstyle\circ}{e}}
\newcommand{\emarui}{\stackrel{\scriptscriptstyle\circ}{e}_i}
\newcommand{\emaruj}{\stackrel{\scriptscriptstyle\circ}{e}_j}
\newcommand{\emaruk}{\stackrel{\scriptscriptstyle\circ}{e}_k}
\newcommand{\emarul}{\stackrel{\scriptscriptstyle\circ}{e}_l}
\newcommand{\ema}{{\stackrel{\scriptscriptstyle\circ}{e}}}
\newcommand{\emai}{{\stackrel{\scriptscriptstyle\circ}{e}_i}}
\newcommand{\emaj}{{\stackrel{\scriptscriptstyle\circ}{e}_j}}
\newcommand{\emak}{{\stackrel{\scriptscriptstyle\circ}{e}_k}}
\newcommand{\emal}{{\stackrel{\scriptscriptstyle\circ}{e}_l}}
\newcommand{\wsta}{{\stackrel{\scriptscriptstyle *}{w}}}
\newcommand{\ehoi}{e_{{\HH}^{\circ},i}}
\newcommand{\ehoj}{e_{{\HH}^{\circ},j}}
\newcommand{\ehok}{e_{{\HH}^{\circ},k}}
\newcommand{\ehol}{e_{{\HH}^{\circ},l}}
\newcommand{\emahoi}{{\stackrel{\scriptscriptstyle\circ}{e}_{{\HH}^{\circ},i}}}
\newcommand{\emahoj}{{\stackrel{\scriptscriptstyle\circ}{e}_{{\HH}^{\circ},j}}}
\newcommand{\emahok}{{\stackrel{\scriptscriptstyle\circ}{e}_{{\HH}^{\circ},k}}}
\newcommand{\emahol}{{\stackrel{\scriptscriptstyle\circ}{e}_{{\HH}^{\circ},l}}}
\newcommand{\epma}{{\stackrel{\scriptscriptstyle\circ}{\varepsilon}}}
\newcommand{\epmai}{{\stackrel{\scriptscriptstyle\circ}{\varepsilon}_i}}
\newcommand{\epmaj}{{\stackrel{\scriptscriptstyle\circ}{\varepsilon}_j}}
\newcommand{\epmak}{{\stackrel{\scriptscriptstyle\circ}{\varepsilon}_k}}
\newcommand{\epmal}{{\stackrel{\scriptscriptstyle\circ}{\varepsilon}_l}}
\newcommand{\End}{\mathrm{End}}
\newcommand{\Aut}{\mathrm{Aut}}
\newcommand{\vEndv}{{}_{\V}\mathrm{End}_{\V}}
\newcommand{\vAutv}{{}_{\V}\mathrm{Aut}_{\V}}
\newcommand{\GLF}{\mathrm{GLF}}
\newcommand{\spa}{\mathrm{span}}
\newcommand{\op}{\mathrm{op}}
\newcommand{\cop}{\mathrm{cop}}
\newcommand{\bop}{\mathrm{bop}}
\newcommand{\Fun}{\mathrm{Fun}}
\newcommand{\HH}{\frak H}
\newcommand{\KK}{\frak K}
\newcommand{\Hhat}{\hat{\frak H}}
\newcommand{\Ht}{\tilde{\frak H}}
\newcommand{\HG}{\frak{H} (\G)}
\newcommand{\Hc}{\mathrm{Hc}}
\newcommand{\HcH}{\mathrm{Hc}({\frak H})}
\newcommand{\Aw}{{\frak A}(w)}
\newcommand{\SS}{{\frak S}(A_{N-1};t)_{\epsilon}}
\newcommand{\SSt}{{\frak S}(A_{N-1};t,\theta)}
\newcommand{\Ss}{{\frak S}(A_{1};t)}
\newcommand{\AS}{{\frak A}(w_{N,t,\epsilon})}
\newcommand{\Str}[1]{\mathrm{Str}^{#1} (\G, *)}
\newcommand{\Rp}{{\cal R}^+}
\newcommand{\Rm}{{\cal R}^-}
\newcommand{\Rpm}{{\cal R}^{\pm}}
\newcommand{\Rmp}{{\cal R}^{\mp}}
\newcommand{\Rph}{{\hat{\cal R}}^+}
\newcommand{\Rmh}{{\hat{\cal R}}^-}
\newcommand{\Rpmh}{{\hat{\cal R}}^{\pm}}
\newcommand{\Rmph}{{\hat{\cal R}}^{\mp}}
\newcommand{\Rt}{{\tilde{\cal R}}}
\newcommand{\Rpt}{{\tilde{\cal R}}^+}
\newcommand{\Rmt}{{\tilde{\cal R}}^-}
\newcommand{\Rpmt}{{\tilde{\cal R}}^{\pm}}
\newcommand{\Qp}{{\cal Q}^+}
\newcommand{\Qm}{{\cal Q}^-}
\newcommand{\Qpm}{{\cal Q}^{\pm}}
\newcommand{\Qmp}{{\cal Q}^{\mp}}
\newcommand{\G}{\EuScript G}
\newcommand{\V}{\EuScript V}
\newcommand{\C}{\Bbb C}
\newcommand{\K}{\Bbb K}
\newcommand{\Z}{\Bbb Z}
\newcommand{\LD}{\mathrm{LD}}
\newcommand{\GLD}{\G_{\mathrm{LD}}}
\newcommand{\wLD}{w_{\mathrm{LD}}}
\newcommand{\GLEH}{\mathrm{GLE} (\HH)}
\newcommand{\aaa}{\bold a}
\newcommand{\bbb}{\bold b}
\newcommand{\ccc}{\bold c}
\newcommand{\ddd}{\bold d}
\newcommand{\p}{\bold p}
\newcommand{\q}{\bold q}
\newcommand{\r}{\bold r}
\newcommand{\s}{\bold s}
\newcommand{\taaa}{\tilde{\bold a}}
\newcommand{\tbbb}{\tilde{\bold b}}
\newcommand{\tccc}{\tilde{\bold c}}
\newcommand{\tddd}{\tilde{\bold d}}
\newcommand{\tp}{\tilde{\bold p}}
\newcommand{\tq}{\tilde{\bold q}}
\newcommand{\tr}{\tilde{\bold r}}
\newcommand{\ts}{\tilde{\bold s}}
\newcommand{\st}{\frak {s}}
\newcommand{\en}{\frak {r}}
\newcommand{\suma}{\sum_{(a)}}
\newcommand{\sumb}{\sum_{(b)}}
\newcommand{\sumc}{\sum_{(c)}}
\newcommand{\sumd}{\sum_{(d)}}
\newcommand{\sumx}{\sum_{(x)}}
\newcommand{\sumy}{\sum_{(y)}}
\newcommand{\sumab}{\sum_{(a),(b)}}
\newcommand{\sumbc}{\sum_{(b),(c)}}
\newcommand{\sumad}{\sum_{(a),(d)}}
\newcommand{\sumxy}{\sum_{(x),(y)}}
\newcommand{\sumax}{\sum_{(a),(x)}}
\newcommand{\sumbx}{\sum_{(b),(x)}}
\newcommand{\sumabc}{\sum_{(a),(b),(c)}}
\newcommand{\sumabcd}{\sum_{(a),(b),(c),(d)}}
\newcommand{\sumk}{\sum_{k \in \V}}
\newcommand{\suml}{\sum_{l \in \V}}
\newcommand{\sumkl}{\sum_{k,l \in \V}}
\newcommand{\face}[4]{\left( \scriptstyle{#1} 
                      \textstyle{\frac[0pt]{#2}{#3}} \scriptstyle{#4} \right)}
\newcommand{\Wrpqs}[4]{w \! \left[ #1 \, \frac[0pt]{#2}{#3} \, #4 \right]}
\newcommand{\wrpqs}[4]{w \! \left[ \scriptstyle{#1} 
                      \textstyle{\frac[0pt]{#2}{#3}} \scriptstyle{#4} \right]}
\newcommand{\hijk}[4]{ \!\! \left[ {#1 \,\, #2} \atop {#3 \,\, #4} \right]}

\title{Coribbon Hopf (face) algebras generated by lattice models}
\author{Takahiro Hayashi}
\date{$\qquad$ Department of Mathematics, 
Furo-cho, Chikusa-ku, Nagoya 464, Japan}
\maketitle

\begin{abstract}
 By studying ``points of the underlying quantum groups''of
 coquasitriangular Hopf (face) algebras, 
 we construct ribbon categories
 for each lattice models without spectral parameter
 of both vertex and face type.
 Also, we give a classification of the braiding and 
 the ribbon structure on quantized classical groups 
 and modular tensor categories closely related to quantum 
 $SU(N)_L$-invariants of 3-manifolds.
\end{abstract}

\section{introduction}

It is widely accepted that (co-)quasitriangular Hopf algebra
is a good algebraic notion which expresses ``quantum groups.''
For, example, each lattice model $w$ of vertex type (and of
face type) without spectral parameter naturally generates
a coquasitriangular (CQT) Hopf (face) algebra, 
thanks to the FRT construction and the Hopf closure (or Hopf
envelope) construction.
The former construction assigns $w$ to the CQT bialgebra (or face algebra)
$\Aw$ (cf. \cite{RTF}, \cite{LarsonTowber}, \cite{Schauenburg}, \cite{gd}),
while the latter construction assigns some CQT bialgebra (or face algebra) $\HH$
to the CQT Hopf (face) algebra $\Hc (\HH)$ (\cite{Phung}, \cite{gsg}). 
However, to give applications of CQT Hopf (face) algebras $\HH$ to 
low-dimensional topology, we need one additional structure on these, 
which is called the {\it ribbon functional} on $\HH$, a dual notion of
the ribbon element.
It is known that there exists a Drinfeld's double of a finite-dimensional 
Hopf algebra, which has no ribbon element (cf. \cite{KauffmanRadford} 
Proposition 7).
Also, it is known that the ribbon functional of a CQT Hopf 
algebra is not necessarily unique even if it exists.
Hence it is natural to investigate sufficient conditions 
for the existence of the ribbon functional
on CQT Hopf (face) algebras, and to develop the classification theory
of the ribbon functionals.

One of the purpose of this paper is to prove the existence 
of the ribbon functionals on 
CQT Hopf face algebras of the form $\Hc (\Aw)$
(cf. Theorem \ref{existrib}). 
This result implies that each $w$ produces a ribbon category,
and therefore, it implies that $w$ generates a family of link invariants.  
We note that when $w = \Check{R}$ is a vertex model, 
this family contains the link invariant constructed by
Reshetikhin \cite{Reshetikhin} (see Remark at the end of Section 6). 
As byproducts, we also obtain several useful results on
the ribbon functionals on CQT Hopf face algebras.

The other purpose of this paper is to
give the classification of the braidings and 
the ribbon functionals on the function algebras
$\Fun (G_q)$ of the quantized classical groups 
$G_q = GL_q(N)$, $SL_q(N)$, $SO_q(N)$, 
$O_q(N)$ and $Sp_q(N)$,
and also, on some Hopf face algebras $\SS$
which are closely related to the $SU(N)_L$-topological
quantum field theories.				  
The braiding of these Hopf (face) algebras is not unique.
However, the non-uniqueness is explained using certain 
gradings of these algebras via cyclic groups $\Gamma$. 
The ribbon functionals of these algebras always exist
and the number of those is at most two.
We note that the proof of the former is very similar to
that of the classification of the braidings of
$\mathrm{Fun} (\mathrm{Mat}_q (N))$ due to Takeuchi
\cite{Takeuchicocycle}, 
while the proof of the latter essentially depends on
our general theory for $\Hc (\Aw)$.   

The paper is organized as follows.
In Section 2, we recall basic concepts of the face algebra.
The notion of the face algebra generalizes that of 
the bialgebra, and is necessary to study lattice models of
face type and the corresponding link invariants 
in the framework of the quantum group theory. 
In Section 3, we recall the relation between lattice models   
and face algebras.
In Section 4, we recall the Hopf closure construction
which is the the main tool of this paper.
In Section 5, we give a study of group-like elements of 
the dual algebras $\Aw^{\circ}$ and $\Hc (\HH)^{\circ}$.
It plays a crucial role for our study of the ribbon functionals.
In Section 6, we give several results on the ribbon functionals
of $\Hc (\HH)$ and its quotients.
In Section 7 and Section 8, we give the results for 
quantized classical groups and the algebras $\SS$ stated above.

The author is grateful to Dr. S. Suzuki for his useful informations.





%
\begin{equation*}
\end{equation*}
\par\noindent
%
{\it Notation}.
Throughout this paper, we use Sweedler's sigma notation for 
coalgebras $C$ and their right comodules $U$, such as 
$(\Delta \otimes \mathrm{id}) (\Delta (a)) =
\sum_{(a)} a_{(1)} \otimes a_{(2)} \otimes a_{(3)}$
(cf. \cite{Montgomery}).
Also, we denote by 
$\rho_U$ the coaction $U \to U \otimes C;$ 
$u \mapsto \sum_{(u)} u_{(0)} \otimes u_{(1)}$,
and by $\pi_U$ 
the left action of $C^*$ on $U$ given by
$\pi_U (X) u = \sum_{(u)} u_{(0)} X(u_{(1)})$
$(u \in U, X \in C^*)$.
For a linear operator $A$ on a vector space $W$ with basis
$\{ \p \}$, we define its matrix $[A^{\p}_{\q}]_{\p \q}$ by 
$A \q = \sum_{\p} \p A^{\p}_{\q}$. 
%

\section{face algebras}

						
Let $\HH$ be an algebra over a field $\K$ equipped
with a coalgebra structure $(\HH,\Delta,\cu)$.  
Let ${\EuScript V}$ be a finite nonempty set and 
let $e_{\HH,i} = e_i$ and ${\ema_{\HH,i} = \emai}$ 
$(i \in {\V})$ be elements of ${\HH}$. 
We say that $({\HH}, \{ e_i,\emai \})$ is a {\em ${\V}$-face algebra} 
if the following relations are satisfied:
\begin{equation}
 \Delta (ab) = \Delta(a) \Delta(b),
\label{D(ab)} 
\end{equation}
\begin{equation}
 e_ie_j = \delta_{ij}e_i, 
 \quad \emai \emaj = \delta_{ij} \emai,
 \quad e_i \emaj = \emaj e_i,
\label{ee} 
\end{equation}
\begin{equation}
 \sum_{k \in \V} \, e_k = 1 = \sum_{k \in {\V}} {\emaruk}, 
\label{sume} 
\end{equation}
\begin{equation}
 \Delta(\emai e_j) = \sum_{k \in {\V}} \emai e_k \otimes \emak e_j,  
\quad
  \cu (\emai e_j) = \delta _{ij} ,  
\label{D(ee)}
\end{equation}
\begin{equation}
 \cu (ab) = \sum_{k \in {\V}} \cu (ae_k) \cu (\emak b)
\label{e(ab)}
\end{equation}
for each $a,b \in \HH$ and $i,j \in \V$.
We call elements $e_i$ and $\emai$ {\it face idempotents}
of $\HH$.
We denote by $\frak{E} = \frak{E}_{\HH}$
the subalgebra of $\HH$ generated by face idempotents.
It is known that bialgebra is an equivalent notion of 
$\V$-face algebra with $\#(\V) = 1$.
For a $\V$-face algebra, we have the following formulas:
\begin{equation}
\label{e(eae)}
 \cu (\emai a) = \cu (e_i a), 
 \quad \cu (a \emai) = \cu (a e_i), 
\end{equation}
\begin{equation}
\label{ae(eae)}
 \suma a_{(1)} \cu  (e_i a_{(2)} e_j) = e_i a e_j,  
\end{equation} 
\begin{equation}
\label{e(eae)a}
 \suma \cu (e_i a_{(1)} e_j)a_{(2)} = \emai a \emaj,  
\end{equation}
\begin{equation}
\label{D(a)}
 \Delta (a) = \sumkl \suma 
e_k a_{(1)} e_l \otimes \emak a_{(2)} \emal,  
\end{equation}
\begin{equation}
\label{eae*a}
 \suma e_i a_{(1)} e_j \otimes a_{(2)} = 
 \suma a_{(1)} \otimes \emai a_{(2)} \emaj,  
\end{equation}
\begin{equation}
\label{D(eeaee)}
 \Delta (\emai e_j a \ema_{i'} e_{j'}) =
 \suma \emai a_{(1)} \ema_{i'} \otimes e_j a_{(2)} e_{j'}  
\end{equation}
for each $a \in \HH$ and $i,j,i',j' \in \V$.

%
%
For a  ${\V}$-face algebra $\HH$, 
its {\it dual face algebra} ${\HH}^{\circ}$
\cite{fdd} is defined to be the dual coalgebra of $\HH$
equipped with product and face idempotents 
given by $\langle XY,\, a \rangle$ $=$
$\sum_{(a)} \langle X,\, a_{(1)} \rangle \langle Y,\, a_{(2)} \rangle$
$(X, Y \in \HH^{\circ}, a \in \HH)$ and
\begin{equation}
\label{eHoi}
 \langle e_{{\HH}^{\circ},i},\, a \rangle
 = \cu (a e_{\HH,i}),
 \quad				
 \langle \ema_{{\HH}^{\circ},i},\, a \rangle
 = \cu (e_{\HH,i} a)
 \quad (a \in \HH, i \in \V).
\end{equation}

Let $x^+$, $x^-$, $e^+$ and $e^-$ be elements of an arbitrary algebra $A$.
We say that $x^-$ is an {\em $(e^+,e^-)$-generalized inverse} of $x^+$
if the following four relations are satisfied:

\begin{equation}
x^{\mp} x^{\pm} = e^{\pm}, 
\quad x^{\pm} x^{\mp} x^{\pm} = x^{\pm}. 
\label{gen.inv.def}
\end{equation}
We note that the $(e^+,e^-)$-generalized inverse of $x^+$
is unique if it exists.

We say that a linear map $S\!: {\HH} \to {\HH}$ is an {\em antipode} of ${\HH}$, 
or ${\HH}$ is a  {\em Hopf} ${\V}$-{\em face algebra} if 
$S$ is the $(E^+, E^-)$-generalized inverse of $\mathrm{id}_{\HH}$
with respect to the convolution product of $\End_{\K} (\HH)$, where 
\begin{equation}
 \label{E+adef}
 E^+ (a) = \sumk \cu (ae_k) e_k, 
  \quad				
 E^- (a) = \sumk \cu (e_k a) \emak
 \quad (a \in \HH).
\end{equation}
%
%
%
%
%
An antipode of a $\V$-face algebra is 
an antialgebra-anticoalgebra map, which satisfies
\begin{equation}
 S(\emai e_j) = \emaj e_i
 \quad (i,j \in {\V}).
\label{S(ee)}
\end{equation}
Let ${\HH}$ be a ${\V}$-face algebra and 
let $\Rp = \Rp_{\HH}$ be an element of $({\HH} \otimes{\HH})^*$ with 
$(m^* (1), (m^{\text{op}})^* (1))$-generalized inverse $\Rm = \Rm_{\HH}$, 
where $m\!: {\HH} \otimes {\HH} \to {\HH}$ denotes the product of ${\HH}$.
We say that 
$\Rp$ is a  {\em braiding} of $\HH$ or 
$(\HH,\Rpm)$ is a {\em coquasitriangular}
(or {\em CQT}) ${\V}$-face algebra 
if the following relations are satisfied:
\begin{equation}
\label{RmXR}
 \quad
 \Rp m^*(X) \Rm = (m^{\text{op}})^*(X)
 \quad (X \in \HH^*),
\end{equation}
\begin{equation}
\label{mid(R)}
 (m {\otimes}\text{id})^* (\Rp) = \Rp_{13} \Rp_{23}, 
 \quad
 (\text{id}{\otimes}m)^* (\Rp) = \Rp_{13} \Rp_{12}.   
\end{equation}
Here for $Z \in (\HH \otimes \HH)^*$ and
$\{ i,j,k \} = \{1,2,3 \}$, we define 
$Z_{ij} \in (\HH^{\otimes 3})^*$ by
$Z_{ij} (a_1, a_2, a_3) = Z(a_i, a_j) \cu (a_k)$  
$(a_1, a_2, a_3 \in \HH)$. 
The braiding $\Rp$ satisfies 
the following relations:
\begin{equation}
 \Rpm_{12} \Rpm_{13} \Rpm_{23} = \Rpm_{23} \Rpm_{13} \Rpm_{12},  
\label{RRR}
\end{equation} 
\begin{equation}
 \cal{R}^{\mp}_{23} \Rpm_{12} \Rpm_{13} 
 = \Rpm_{13} \Rpm_{12} \cal{R}^{\mp}_{23}, 
\quad
 \cal{R}^{\mp}_{13} \cal{R}^{\mp}_{23} \Rpm_{12} 
 = \Rpm_{12} \cal{R}^{\mp}_{23} \cal{R}^{\mp}_{13}, 
\label{RmpRR}
\end{equation} 
\begin{equation}
 \Rp (\emai e_j a \emak e_l,\, b) = \Rp (a,\, \emaj e_l b \emai e_k),
\label{Rpeeee}						
\end{equation}
\begin{equation}
 \Rm (\emai e_j a \emak e_l,\, b) = \Rm (a,\, \emak e_i b \emal e_j),
\label{Rmeeee}
\end{equation}
\begin{equation}
 \Rp (\emai e_j,\, a) = \cu (e_j a e_i),
 \quad 
 \Rp (a,\, \emai e_j) = \cu (e_i a e_j),
 \label{Rpeea}
\end{equation}
\begin{equation}
 \Rm (\emai e_j,\, a) = \cu(e_i a e_j),
 \quad 
 \Rm (a,\, \emai e_j) = \cu(e_j a e_i)
\label{Rmeea}
\end{equation}
for each $a,b \in \HH$ and $i,j,k,l \in \V$. 
If $\HH$ is a Hopf $\V$-face algebra, then we have:
\begin{equation}
 (S{\otimes}\text{id})^* (\Rp) = \Rm,   
 \quad 
 (\text{id}{\otimes}S)^* (\Rm) = \Rp.
\label{Sid(R)}
\end{equation}
\begin{prop}
\label{altbraiding}
 Let $(\HH, \Rpm)$ be a CQT $\V$-face algebra.\\
 \rom{(1)}
 Then,  $\Rm_{21}\!: a \otimes b \mapsto \Rm(b, a)$ gives
 another braiding of $\HH$.\\
 \rom{(2)}
 Let $\Gamma$ be a semigroup and 
 $\HH = \bigoplus_{\gamma \in \Gamma} \HH_{\gamma}$ 
 a decomposition of $\HH$ such that 
 $\HH_{\gamma} \HH_{\delta} \subset \HH_{\gamma \delta}$ 
 and that 
 $\Delta ( \HH_{\gamma} ) \subset \HH_{\gamma} \otimes \HH_{\gamma}$.
 Let $\chi\!: \Gamma \times \Gamma \to \K^{\times}$ 
 be a map such that
 $\chi (\gamma_1 \gamma_2, \delta)$ $=$		
 $\chi (\gamma_1, \delta) \chi (\gamma_2, \delta)$,
 $\chi (\gamma, \delta_1 \delta_2)$ $=$
 $\chi (\gamma, \delta_1) \chi (\gamma, \delta_2)$.
 Then, there exists a new braiding $\Rp_{\chi}$ of $\HH$ given by 
 \begin{equation}
 \label{Rchidef} 
  \Rpm_{\chi} (a, b) =
  \chi (\gamma, \delta)^{\pm 1} 
  \Rpm (a, b) 
  \quad (a \in \HH_{\gamma}, b \in \HH_{\delta}).
 \end{equation}		
 If $(\HH, \Rpm)$ is closable, then so is $(\HH, \Rpm_{\chi})$. 
\end{prop}
\begin{pf}
 This is straightforward.
\end{pf}
%
%
%
%
Let $({\HH},\Rpm)$ be a CQT Hopf $\V$-face algebra
and $\cal{V}$ an invertible central element of 
${\HH}^*$.
We say that $\cal{V}$ is a {\em ribbon functional} of $\HH$, 
or $(\HH,\cal{V})$ is a
{\em coribbon Hopf} $\V$-{\em face algebra}
if 
%
%
\begin{equation}
 \label{m(V)} 
  m^*({\cal V}) = \Rm \Rm_{21} ({\cal V} \otimes {\cal V}),
\end{equation} 
\begin{equation}
\label{S(V)}
  S^* (\cal{V}) = \cal{V}.
\end{equation}

A map $f\!: \HH \to \KK$ between $\V$-face algebras is called a
{\em map of $\V$-face algebras} if it is both an algebra and 
a coalgebra map such that $f(e_i) = e_i$, $f(\emai) = \emai$
for each $i \in \V$. 
If both $\HH$ and $\KK$ have antipode, then we have
\begin{equation}
 f(S(a)) = S(f(a))
 \quad (a \in \HH).
\label{f(S(a))}
\end{equation} 

A map $f\!: \HH \to \KK$ of $\V$-face algebras 
between CQT $\V$-face algebras is called a
{\em map of CQT $\V$-face algebras} if 
\begin{equation}
 \label{mapCQTdef}
 (f \otimes f)^* (\Rp_{\KK}) = \Rp_{\HH}.
\end{equation}
An ideal $\frak{I}$ of a $\V$-face algebra $\HH$ is called 
a {\it biideal} if it is a coideal of the underlying 
coalgebra of $\HH$.
If in addition, $\HH$ is a CQT $\V$-face algebra and 
$\frak{I}$ satisfies 
$\Rpm(\frak{I},\HH) = \Rpm(\HH, \frak{I}) = 0$,
then $\frak{I}$ is called a {\it CQT biideal} of $\HH$. 
For each $\V$-face algebra (resp. CQT $\V$-face algebra)
$\HH$ and its biideal (resp. CQT biideal) 
$\frak{I}$, the quotient $\HH / \frak{I}$ becomes a 
$\V$-face algebra (resp. CQT $\V$-face algebra)
in an obvious manner.

\section{lattice models and comodules}

Let $\G$ be a finite oriented graph 
with set of vertices ${\V}$ = ${{\G}^0}$.
For an edge ${\p}$, we denote by ${\st (\p)}$ and ${\en (\p)}$
its {\em source} ({\em start}) and its {\em range} ({\em end}) respectively.
For each $m \geq 1$, we denote by 
${{\G}^m = {\coprod}_{i,j \in \V}{\G}_{ij}^m}$ 
the set of {\em paths} of $\G$ of {\em length} $m$, that is, 
$\p \in {\G}_{ij}^m$ if $\p$ is a sequence
$(\p_1, \ldots, \p_m)$ of edges of $\G$ such that
$\st (\p):= \st (\p_1) = i$,
$\en (\p_n) = \st(\p_{n+1})\,$ $(1 \leq n < m)$ 
and $\en (\p):= \en(\p_m) = j$.
Also, we set $\st (i) = i = \en (i)$, $\G_{ii}^0 = \{ i\}$
and $\G_{ij}^0 = \emptyset$
for each $i \in \V$ and $j \ne i$.
Let ${\HH (\G)}$ be the linear span of 
the symbols ${e{\p \choose \q}}$
\( ( \p,\q \in {\G}^m, m \geq 0) \).
Then ${\HH (\G)}$ becomes a $\V$-face algebra by setting
\begin{equation}
 \emai = \sum_{j \in \V}e{i \choose j}, 
 \quad e_j = \sum_{i \in \V}e{i \choose j},  
\label{eHGi}
\end{equation}
\begin{equation}
 e{\p \choose \q} e{\r \choose \s} = 
 {\delta}_{\en(\p) \st(\r)} \, {\delta}_{\en(\q) \st(\s)} 
 \: e{\p\cdot\r \choose \q\cdot\s}, 
\label{epqers}
\end{equation}
\begin{equation}
 \Delta \left( e{\p \choose \q} \right) = 
 \sum_{{\bold t} \in {\EuScript G}^m} 
 e{\p \choose {\bold t}} \otimes e{{\bold t} \choose \q},  
\quad
 \cu \left( e{\p \choose \q} \right)  = {\delta}_{\p \q}
\label{D(epq)}
\end{equation}
for each $\p,\q \in {\G}^m$ 
and $ \r,\s \in {\G}^n$ $( m,n \geq 0)$. 
Here for paths $\p = (\p_1,\ldots,\p_m)$ and 
$\r = (\r_1,\ldots,\r_n)$, we set 
$\p \cdot \r = (\p_1,\ldots,\p_m,\r_1,\ldots,\r_n)$ 
if $\en(\p) = \st(\r)$ and $m, n \geq 1$, 
and also, we set 
$\st (\p) \cdot \p = \p = \p \cdot \en (\p)$
for each $\p \in \G^m$ $(m \geq 0)$.

We say that a quadruple 
$\left( \r \frac[0pt]{\p}{\q} \s \right)$ 
is a {\it face} if 
$\p,\q, \r, \s \in {\G}^1$ and
\begin{equation} 
 \st (\p) = \st (\r),
 \quad \en (\p) = \st (\s),  
 \quad \en (\r) = \st (\q), 
 \quad \en (\q) = \en (\s).
 \label{facecond}
\end{equation}
We say that $(\G, w)$ is a {\it face model} 
(or $\V$-{\it face model} ) over ${\Bbb K}$ if 
$w$ is a map which assigns a scalar
$w \!\! \left[ \r \frac[0pt]{\p}{\q} \s \right] \in {\Bbb K}$
to each face 
$\left( \r \frac[0pt]{\p}{\q} \s \right)$
of $\G$.
A face model $(\G, w)$ is called a {\em vertex model} if 
$\#(\V) = 1$.
For convenience, we set 
$w \!\! \left[ \r \frac[0pt]{\p}{\q} \s \right] = 0$
unless
$\p, \q, \r, \s \in \G^1$ satisfy \eqref{facecond}.
For a face model $(\G,w)$, we identify $w$ 
with the linear operator on
$\K \G^2 := \bigoplus_{\p \in \G^2} \K \p$ given by
\begin{equation}
\label{w(pq)}
 w(\p,  \q) =
 \sum_{(\r, \s) \in {\G}^2}
 w \!\! \left[ \r \frac[0pt]{\p}{\s} \q \right]
 (\r, \s) 
 \quad ( (\p, \q) \in \G^2).
\end{equation}
A face model is called {\it invertible}
if $w$ is invertible as an operator on  
$\K \G^2$.
For an invertible face model $(\G, w)$, 
we define another face model $(\G, w^{-1})$, 
using the identification \eqref{w(pq)}.
An invertible face model is called {\it star-triangular}
(or {\it Yang-Baxter}) if $w$ satisfies
the braid relation $w_1 w_2 w_1 = w_2 w_1 w_2$, 
where $w_1$ and $w_2$ denote linear operators on
$\K \G^3$ defined by
$w_1(\p, \q, \r) =   w(\p, \q) \otimes \r$
and
$w_2(\p, \q, \r) =   \p \otimes w(\q, \r)$.
Here we identify $(\p, \q, \r) \in \G^3$
with 
$\p \otimes \q \otimes \r \in 
(\K \G^1)^{\otimes 3}$.

For a face model $(\G, w)$, we define the algebra 
$\frak A (\G, w)= \Aw$ 
to be the quotient of $\HH (\G)$ modulo the following relations:
\begin{equation}
 \sum_{(\ccc, \ddd) \in {\G}^2}
 w \!\! \left[ \aaa \frac[0pt]{\ccc}{\bbb} \ddd \right] 
 e{\ccc \cdot \ddd  \choose \p \cdot \q} =
 \sum_{(\r, \s) \in {\G}^2}
 w \!\! \left[ \r \frac[0pt]{\p}{\s} \q \right]
 e{\aaa \cdot \bbb  \choose \r \cdot \s}
 \quad ( (\p, \q),\, (\aaa, \bbb) \in {\G}^2).
\label{relAw}
\end{equation}
Then $\Aw$ has a unique structure of 
$\V$-face algebra such that the projection
$\HH (\G) \to \Aw$ is a map of $\V$-face algebras.
For each $n \geq 0$, 
$\K \G^n$ becomes a comodule of  
$\frak{A}_n (w)$
via $\rho (\q) = \sum_{\p \in \G^n} \p \otimes e\!\binom{\p}{\q}$,
where the subcoalgebra $\frak{A}_n (w)$ of $\Aw$ is defined as
the linear span of the elements of the form
$e\!\binom{\p}{\q}$ $(\p, \q \in \G^n)$. 
If $(\G,w)$ is star-triangular, then there exist unique 
bilinear pairings $\Rpm$ on $\Aw$
such that $(\Aw, \Rpm)$ is a CQT $\V$-face algebra
and that 
\begin{equation}
  \Rp \left( e{\p \choose \q},\; e{\r \choose \s} \right)
 =  w \!\! \left[ \r \frac[0pt]{\q}{\p} \s \right]
\label{UnivR=W}
\end{equation}
for each $\p,\q, \r, \s \in {\G}^1$
(cf. \cite{LarsonTowber}, \cite{gd}, \cite{Schauenburg}, \cite{fb}). 
We call $\Rp$ the {\it canonical braiding} of $\Aw$.
For a vertex model $w = \Check{R}$, $\Aw$ coincides with
FRT bialgebra $A_R$, where $R = P \Check{R}$ 
and $P(\p, \q) = (\q, \p)$.

Let $\tilde{\G}$ be the orientation-reversed graph of $\G$
and let $\tilde{}: \G^m \to \tilde{\G}^m \, ;\; \p \mapsto \tilde{\p}$ 
$(m \geq 0)$ be the canonical bijection which satisfies
$\widetilde{\p \cdot \q} = \tilde{\q} \cdot \tilde{\p}$,
$\st(\tilde{\p}) = \en({\p})$ and
$\en(\tilde{\p}) = \st({\p})$.
We also define a new graph $\GLD$
by setting $\GLD^0 = \V$ and
$\GLD^1 = \G^1 \coprod \tilde{\G}^1$.
Let $\G \bar{\times} \tilde{\G}$
and $\tilde{\G} \bar{\times} \G$ 
denote subsets of $\GLD^2$
consisting of elements of the form
$\p \cdot \tilde{\q}$
and $\tilde{\p} \cdot \q$
$(\p,\q \in \G^1)$ 
respectively. 
We define linear operators 
$\wLD, \wLD^-\!: 
\K (\tilde{\G} \bar{\times} \G)
\to \K (\G \bar{\times} \tilde{\G})$
by
\begin{equation}
 \wLD (\tilde{\p} \cdot \q) =
 \sum_{\r,\s} \wLD \!\! \left[ \r \frac[0pt]{\tilde{\p}}{\tilde{\s}} \q \right] 
 \r \cdot \tilde{\s}\, ; 
 \quad
 \wLD \!\! \left[ \r \frac[0pt]{\tilde{\p}}{\tilde{\s}} \q \right] 
 =:
 w^{-1} \!\! \left[ \p \frac[0pt]{\q}{\r} \s \right],  
\end{equation}
\begin{equation}
 \wLD^- (\tilde{\p} \cdot \q) =
 \sum_{\r,\s} \wLD^- \!\! \left[ \r \frac[0pt]{\tilde{\p}}{\tilde{\s}} \q \right] 
 \r \cdot \tilde{\s} \, ; 
 \quad
 \wLD^- \!\! \left[ \r \frac[0pt]{\tilde{\p}}{\tilde{\s}} \q \right] 
 =:
 w \!\! \left[ \p \frac[0pt]{\q}{\r} \s \right]. 
\end{equation}
We say that a star-triangular $\V$-face model $(\G,w)$
is {\it closable}
if both $\wLD$ and $\wLD^-$ are invertible. 
In this case, we define a new $\V$-face model $(\GLD,\wLD)$
by extending $\wLD$ on $\K \GLD^2$ via
$\wLD |_{\K \G^2} = w$,  
$\wLD |_{\K (\G \bar{\times} \tilde{\G})} 
= (\wLD^-)^{-1}$  
and 
\begin{equation}
 \wLD \!\! \left[ \tilde{\r} \frac[0pt]{\tilde{\p}}{\tilde{\q}} \tilde{\s} \right] 
 =
 w \!\! \left[ \s \frac[0pt]{\q}{\p} \r \right]
 \quad (\p, \q, \r, \s \in \G^1).
\end{equation}
We call $\wLD$ the {\it Lyubashenko double} of $w$.
As in case $(\G,w)$ is a vertex model, 
$(\GLD,\wLD)$ is a star-triangular face model.

Let $\HH$ be a $\V$-face algebra and $U$ its (right) comodule.
We define its {\it face space decomposition}
$U = \bigoplus_{i,j \in {\V}}U(i,j)$ by
$U(i,j) = \pi_U (\ema_i e_j) (U)$.
%
%
Let $V$ be another $\HH$-comodule. We define the 
{\it truncated tensor product} 
$U \overline{\otimes} V$ to be the vector space
\begin{equation}
 U \overline{\otimes} V = {\bigoplus}_{i,j,k \in {\V}}U(i,k) \otimes V(k,j) 
\end{equation}								
equipped with the $\HH$-comodule structure given by
\begin{equation}
 \rho_{U \bar{\otimes} V} (\overline{u \otimes v} )
 = \sum_{(u),(v)} 
 \left( u_{(0)} \otimes v_{(0)} \right)
 \otimes u_{(1)}v_{(1)},
\end{equation}
where $\bar{}\!: U \otimes V \to U \bar{\otimes} V$
denotes the projection 
$\sum_k \pi_U (e_k) \otimes \pi_V (\emak) $.
For $\HH$-comodules $U, U^{\prime},$ $V, V^{\prime}$ and
maps 
$f \in \End_{\pi (\frak{E})} (U, U^{\prime})$,
$g \in \End_{\pi (\frak{E})} (V, V^{\prime})$,
we set				
\begin{equation}
 f \bar{\otimes} g = 
 (f \otimes g) |_{U \overline{\otimes} V},
\end{equation}
where $\frak{E} = \frak{E}_{\HH^{\circ}}$.
If both $f$ and $g$ are comodule maps, then  
so is $f \bar{\otimes} g$.  
The category $\bold{Com}_{\HH}$ of all 
$\HH$-comodules becomes a monoidal category via 
$\bar{\otimes}$ and
the category $\bold{Com}^f_{\HH}$ of all finite-dimensional
$\HH$-comodules becomes its sub monoidal category.
The category $\bold{Com}^f_{\HH}$ is rigid
if and only if 
$\HH$ has a bijective antipode.

Next, suppose $\HH$ has a braiding $\Rpm$.
Then $\bold{Com}_{\HH}$ and $\bold{Com}^f_{\HH}$ 
become braided categories
via the functorial isomorphism 
$c_{UV}\!: U \overline{\otimes} V$ 
$\cong V \overline{\otimes} U$ given by
\begin{equation}
 c_{UV} (\overline{u \otimes v}) 
 = \sum_{(u),(v)} v_{(0)} \otimes u_{(0)} 
 {\cal R}^+(u_{(1)},v_{(1)}).
\end{equation}
If, in addition, $\HH$ has a ribbon functional $\cal{V}$, 
then $\bold{Com}^f_{\HH}$ becomes a ribbon category
(see e.g. \cite{Kassel})
via twist
$\theta_{U}\!: U \cong U$ given by
$\theta_{U} = \pi_U (\cal{V}^{-1})$.
Conversely, we have the following.

\begin{prop}[\cite{LarsonTowber}, \cite{Kassel}]
\label{correspalgcom}
 Let $\HH$ be a $\V$-face algebra such that 
 either $\bold{Com}_{\HH}$ or
 $\bold{Com}^f_{\HH}$ is a braided monoidal category
 with braiding $\{ c_{UV} \}$.
 Then, $\HH$ becomes a CQT $\V$-face algebra via
 \begin{gather}
 \label{Rp=eec}
  \Rp (a, b) = 
  \sumkl ( \cu \otimes \cu ) \circ
  c_{LM} (\emak a e_l \otimes e_l b \emak ), \\
 \label{Rm=eec} 
  \Rm (b, a) = 
  \sumkl ( \cu \otimes \cu ) \circ
  (c_{ML})^{-1} (\emak a e_l \otimes e_l b \emak ), 
 \end{gather}
 where $L$ and $M$ denote arbitrary finite-dimensional
 sub $\HH$-comodules of $\HH$ such that 
 $\emai a e_j \in \L$,  $e_j b \emai \in M$
 $(i, j \in \V)$.
 If, in addition, $\HH$ has an antipode and 
 $\bold{Com}^f_{\HH}$ is a ribbon category
 with twist $\{ \theta_{U} \}$, then $\HH$ becomes a 
 coribbon Hopf $\V$-face algebra via 
 \begin{equation}
 \label{V=etheta-}
  \cal{V}^{\pm 1} (a) = 
  \cu ( \theta_L^{\mp 1} (a) ).
 \end{equation}
\end{prop}
\begin{pf}
 To begin with, we note that the existence of such 
 $L$ and $M$ follows from 
 the fundamental theorem of coalgebras, 
 and that \eqref{Rp=eec}-\eqref{V=etheta-}
 do not depend on the choice of 
 $L$ and $M$ because of the naturality of $c$ and $\theta$.
 Here, we will give a proof of the last assertion.
 Let $\cal{V}^{\pm} \in \HH^*$ be as in  \eqref{V=etheta-}
 and let $U$ be a finite-dimensional $\HH$-comodule.
 For each $u^* \in U^*$, we define the 
 $\HH$-comodule map $F_{u^*}\!:U \to \HH$ by
 $F_{u^*}(u)$ $=$ $\sum_{(u)} \langle u^*, u_{(0)} 
 \rangle u_{(1)}$ $(u \in U)$.
 Then, we have
 \begin{align}
 \label{u*thetam}	
  \langle u^*,\, \theta_U^{\mp 1} (u) \rangle 
  = &
  \cu \circ F_{u^*} \circ \theta_U^{\mp 1} (u)
  =
  \cu \circ \theta_{\mathrm{Im} (F_{u^*})}^{\mp 1} \circ F_{u^*} (u) \\
  = &					
  \langle u^*,\, \pi_U (\cal{V}^{\pm 1}) u \rangle,
  \nonumber
 \end{align}
 or equivalently,
 \begin{equation}
  \label{theta=piV-}
  \theta_U^{\mp 1} (u) = 
  \pi_U (\cal{V}^{\pm 1}) u.
 \end{equation}
 Rewriting 
 $\langle u^*,\, \theta_U (u) \rangle$ $=$
 $\langle \theta_{U^{\lor}} (u^*), u \rangle$
 via this equality, we obtain
 $S (\cal{V}) = \cal{V}$.

 Let $a$ and $b$ elements of $\HH$ and let $L$ 
 and $M$ be as above.
 Since $\overline{a \otimes b} = 
 \sum_k a e_k \otimes e_k b$
 by \eqref{e(ab)}-\eqref{ae(eae)}, we have
 \begin{equation}
 \pi_{L \bar{\otimes} M } (\cal{V}) (\overline{a \otimes b})
 =
 \sum_{(a), (b)} a_{(1)} \otimes b_{(1)}  
 \langle \cal{V}, a_{(2)} b_{(2)}  
 \rangle 
 \end{equation}
 by \eqref{D(eeaee)}.
 Using \eqref{theta=piV-} and the equality 
 $\theta_{L \overline{\otimes} M}^{-1}$ $=$
 $c_{L \overline{\otimes} M}^{-1} \circ
  c_{M \overline{\otimes} L}^{-1} \circ
 (\theta_{L}^{-1} \bar{\otimes} \theta_{M}^{-1})$,
 we see that the left-hand side of the above equality is
 \begin{align}
  \label{cc(ab)VaVb}
  &\sum_{(a), (b)}
  c_{L \overline{\otimes} M}^{-1} \circ
  c_{M \overline{\otimes} L}^{-1} 
  (\overline{a_{(1)} \otimes b_{(1)}}) 
  \cal{V} (a_{(2)}) \cal{V} (b_{(2)}) \\
  = &
  \sum_{(a), (b)} a_{(1)} \otimes b_{(1)}  
  \Rm (a_{(2)}, b_{(2)}) \Rm (b_{(3)}, a_{(3)})  
  \cal{V} (a_{(4)}) \cal{V} (b_{(4)}),
 \end{align}
 where \eqref{cc(ab)VaVb} follows from the fact that
 $\theta_L$ and $\theta_M$ commute with the action of 
 the face idempotents of $\HH^{\circ}$.
 Taking the image via $\cu \otimes \cu$, we get
 \eqref{m(V)}.
\end{pf}

Let $U$ be a finite-dimensional comodule of a CQT
$\V$-face algebra $\HH$. 
For each $i,j \in \V$, choose a basis $\G^1_{ij}$ of $U (i,j)$.
Let $\G$ be the oriented graph with set of vertexes $\V$
and the set of edges $\G^1: = \coprod_{ij} \G^1_{ij}$.
Then we obtain a star-triangular 
$\V$-face model $(\G, w_U)$ be setting
\begin{equation}
\label{wUdef}
 c_{UU} \left( \p \otimes \q \right)
 = \sum_{(\r,\s) \in \G^2} 
 w_U \!\! \left[ \r \frac[0pt]{\p}{\s} \q \right] 
 \r \otimes \s.
 \quad ((\p, \q ) \in \G^2).
\end{equation}
%


\section{Drinfeld functionals and the Hopf closure}

%
%
Let $(\HH,\Rpm)$ be a CQT ${\V}$-face algebra.
We say that ${\HH}$ is {\em closable}
(or ${\HH}$ is a {\em CCQT} ${\V}$-face algebra)
if there exist both 
$(\cal{F}^+, \cal{F}^-)$-generalized inverse $\Qm$ of $\Rp$ and 
$(\cal{F}^-, \cal{F}^+)$-generalized inverse $\Qp$ of $\Rm$
in the algebra $({\HH}\otimes{\HH}^{\cop})^*$, 
where $\cal{F}^{\pm}$ denote bilinear forms on $\HH$
defined by 

\begin{equation}
 \cal{F}^+ (a,\, b) = \sumk \cu (e_k a) \cu(e_k b),
 \quad
 \cal{F}^- (a,\, b) = \sumk \cu (a e_k) \cu(b e_k)
 \quad (a,b \in \HH).
\end{equation}   
We call $\Qpm$  {\em Lyubashenko forms} of $\HH$.
%
%
The Lyubashenko forms of a CQT $\V$-face algebra
are unique if they exist. 
If $\HH$ has an antipode, then $\HH$ is closable with
Lyubashenko forms given by
\begin{equation}
 \label{Qp=Rm}
 \Qp(a,b) = \Rm (S(a), b),
 \quad \Qm(a,b) = \Rp (a, S(b))
 \quad (a,b \in \HH).
\end{equation} 
%
%
%

For a star-triangular face model $(\G,w)$,
$\Aw$ is closable if and only if
$(\G,w)$ is closable.
In this case, Lyubashenko forms
$\Qpm$ of $\Aw$ satisfy
\begin{equation}
  \Qp \left( e{\p \choose \q},\; e{\r \choose \s} \right)
 =  \wLD^{-1} \!\! \left[ \tilde{\q} \, \frac[0pt]{\s}{\r} \, \tilde{\p} \right],
\quad
  \Qm \left( e{\p \choose \q},\; e{\r \choose \s} \right)
 =  \wLD \!\! \left[ \tilde{\s} \, \frac[0pt]{\q}{\p} \, \tilde{\r} \right]
\label{Qp=W}
\end{equation}
for each $\p,\q, \r,\s \in {\G}^1$.

For a CCQT ${\V}$-face algebra ${\HH}$, 
we define linear functionals ${\cal U}_{\nu}$ 
$(\nu = 1,2)$ on ${\HH}$ via 
\begin{equation}
 {\cal U}_1 (a) = \suma \Qm (a_{(2)},a_{(1)}),
\quad
 {\cal U}_2 (a) = \suma \Qp (a_{(1)},a_{(2)})
\quad (a \in \HH)
\label{Udef}
\end{equation}
and call them {\em Drinfeld functionals} of ${\HH}$. 
The Drinfeld functionals of a 
CCQT $\V$-face algebra $\HH$ are invertible
in ${\HH}^*$ and satisfy the following relations\rom{:}
 
\begin{equation}
\label{U-}
 {\cal U}_1^{-1}(a) = \suma \Qp (a_{(2)},a_{(1)}),
\quad
 {\cal U}_2^{-1}(a) = \suma \Qm (a_{(1)},a_{(2)}),
\end{equation}
\begin{equation}
\label{UURp}
 ({\cal U}_{\nu} \otimes {\cal U}_{\nu}) \Rpm = 
 \Rpm ({\cal U}_{\nu} \otimes {\cal U}_{\nu}), 
\end{equation}
\begin{equation}
\label{U(ee)}
 {\cal U}_{\nu}^{\pm 1}({\emarui} e_j) = {\delta}_{ij}, 
\quad
 {\cal U}_{\nu}^{\pm}(\emai a \emaj) =   {\cal U}_{\nu}^{\pm}(e_i a e_j), 
\end{equation} 
\begin{equation}
 \label{U1U2}
 {\cal U}_1 {\cal U}_2 = {\cal U}_2 {\cal U}_1,
\end{equation}
\begin{equation}
\label{m(U)} 
 m^*({\cal U}) = \Rm \Rm_{21} ({\cal U} \otimes {\cal U})
 = ({\cal U} \otimes {\cal U}) \Rm \Rm_{21}, 
\end{equation} 
\begin{equation}
 m^*({\cal U}^{-1}) 
 = \Rp_{21} \Rp ({\cal U} \otimes {\cal U} )^{-1}
 = ({\cal U} \otimes {\cal U} )^{-1} \Rp_{21} \Rp
\label{m(Um)}
\end{equation}
for each $\nu = 1,2$, $a \in \HH$
and $i,j \in \V$, where $\cal{U}$ stands for 
$\cal{U}_1$ or $\cal{U}_2^{-1}$.

Let $f\!: \HH \to \KK$
be a map of CQT $\V$-face algebras.
If $\KK$ is closable with Lyubashenko forms $\Qpm_{\KK}$
and Drinfeld functionals
${\cal{U}_{\nu}}_{\KK}$, 
then $\HH$ is also closable with Lyubashenko forms and
Drinfeld functionals given by
\begin{equation}  
\label{f(U)=U}
 \Qpm_{\HH} = (f \otimes f)^* (\Qpm_{\KK}),
\quad 	
 {\cal{U}_{\nu}}_{\HH} =
 f^* ({\cal{U}_{\nu}}_{\KK}).
\end{equation}			
Next, we recall the {\it Hopf closure} 
(or {\it Hopf envelope}) construction
of CQT Hopf $\V$-face algebras.
It is introduced by Phung Ho Hai \cite{Phung} for 
bialgebras, and independently, by \cite{gsg} for face algebras. 
Let $\HH$ be a CCQT $\V$-face algebra.
We denote by ${\HH}^{\bop}$ its biopposite 
$\V$-face algebra, that is, ${\HH}^{\bop}$ 
is a $\V$-face algebra equipped with
the opposite product and the opposite coproduct
of $\HH$ together with the face idempotents 
\begin{equation}		
 \ema_{\HH^{\bop},i} =  e_{\HH,i},  
\quad
 e_{{\HH}^{\bop},i} =
 \ema_{{\HH},i}.
\end{equation}
Let $\sigma\! : \HH \to {\HH}^{\bop}$ be
the canonical anti-isomorphism, which satisfies
\begin{equation}
 \sigma(e_{\HH,i}) = 
 \ema_{{\HH}^{\bop},i},
\quad
 \sigma(\emaru_{\HH,i}) = 
 e_{{\HH}^{\bop},i}
 \quad (i \in \V).
\end{equation}
Then 
\begin{equation}		
 \hat{\HH} :=
 \HH \otimes_{\frak{E}} \HH^{\bop} =
 \bigoplus_{k,l \in \V}
 \HH \ema_k e_l \otimes \sigma (\HH \emal e_k) 
\end{equation}
becomes a $\V$-face algebra by setting
\begin{equation}
 (a \tensE \sigma(b))(c \tensE \sigma(d)) = 
 \sumbc \Rm (b_{(1)},c_{(3)}) 
 \Qp (b_{(3)},c_{(1)}) a c_{(2)} \tensE \sigma(d b_{(2)}), 
\label{asbcsd}
\end{equation}
\begin{equation}
 \Delta (a \tensE \sigma(b)) = 
 \sumab (a_{(1)} \tensE \sigma(b_{(2)})) 
 \otimes (a_{(2)} \tensE \sigma(b_{(1)})), 
\end{equation}
\begin{equation}
 \varepsilon (a \tensE \sigma(b)) = 
 \sumk {\cu (a e_k) \cu (b e_k)},
\end{equation}
\begin{equation}
 e_{\Hhat,i} = e_{\HH,i} \tensE \sigma(1_{\HH}),
 \quad
 \emaru_{\Hhat,i} = \emaru_{\HH,i} \tensE \sigma(1_{\HH})
\end{equation}
for each $a,b,c,d \in \HH$ and $i \in \V$,
Let ${\frak J}$ be the ideal of $\Hhat$ generated by 
all elements of the form:
\begin{equation} 
 \suma (1 \tensE \sigma(a_{(1)}))(a_{(2)} \tensE 1) 
 - \sumk \varepsilon (a e_k) e_k, 
 \label{Idef1} 
\end{equation}
\begin{equation}
 \suma a_{(1)} \tensE \sigma(a_{(2)}) - \sumk \varepsilon (e_k a) \emak 
 \quad (a \in \HH).
 \label{Idef2} 
\end{equation} 
%
%
%
It is easy to verify that $\frak{J}$ becomes a biideal.
We denote the quotient $\V$-face algebra $\Hhat / {\frak J}$ 
by $\HcH$ and call it the {\em Hopf closure} of $\HH$.
For simplicity, we denote an element $a \tensE \sigma(b) + \frak{J}$ 
of $\HcH$ by $a \sigma(b)$ for each
$a,b \in \HH$.
The Hopf closure $\HcH$ has a unique structure of CQT Hopf 
$\V$-face algebra such that the canonical map
$\iota\!: \HH \to \HcH;$  
$a \mapsto a \tensE 1 + \frak{J}$ $(a \in \HH)$
is a map of CQT $\V$-face algebras.
Explicitly, the antipode of $\HcH$ is given by
\begin{align}			
\label{SHcH}
 S(a \sigma(b)) 
 =  \sumb \cal{U}_{\nu} (b_{(1)}) 
 b_{(2)} \sigma (a) \cal{U}_{\nu}^{-1} (b_{(3)})
 \quad (\nu = 1,2). 
\end{align}
When $\HH$ is a bialgebra, the underlying Hopf algebra of
$\HcH$ agrees with the Hopf envelope of $\HH$ in the sense of 
Manin \cite{Manin}.
%
%
%
The Hopf closure has the following universal mapping property.
\begin{thm}
\label{UMP}					
 Let $\HH$ be a CCQT $\V$-face algebra and $\KK$ 
 a CQT Hopf $\V$-face algebra.
 Let $f\!: \HH \to \KK$ be a map of CQT $\V$-face algebras.
 Then there exists a unique map 
 $\bar{f}\!:\HcH \to \KK$ 
 of CQT $\V$-face algebras
 such that $f = \bar{f} \circ \iota$,
 where $\iota\!: \HH \to \HcH$  
 is given by 
 $\iota (a) = a \tensE 1 + \frak{J}$ $(a \in \HH)$.
 Explicitly, we have 
 \begin{equation}
 \label{barf(asb)}
 \bar{f} (a \sigma (b)) =
 f (a) S(f (b)).
 \end{equation}
\end{thm}
\begin{prop}
\label{univA}
 Let $\HH$ be a CQT $\V$-face algebra
 \rom{(}resp. CQT Hopf $\V$-face algebra\rom{)} 
 and $U$ its finite-dimensional comodule.
 Let $(\G, w_U)$ be a face model given by \eqref{wUdef}.
 Then there exists a unique map $f\!: \frak{A}(w_U) \to \HH$
 \rom{(}resp. $f\!: \Hc(\frak{A}(w_U)) \to \HH$\rom{)}
 of CQT $\V$-face algebras such that
 $(\mathrm{id}_{\K \G^1} \otimes f) \circ \rho_{\HH}$
 $=$ $\rho_{\frak{A}(w_U)}$ \rom{(}resp.
 $(\mathrm{id}_{\K \G^1} \otimes f) \circ \rho_{\HH}$
 $=$ $\rho_{\Hc(\frak{A}(w_U))}$\rom{)}.
\end{prop}
\begin{pf}
 See \cite{fb} for a proof of
 the assertion for $\frak{A}(w_U)$.
 The assertion for $\Hc( \frak{A}(w_U) )$ follows from that of
 $\frak{A}(w_U)$ and the universal mapping property of $\Hc$.
\end{pf}
\begin{prop} 
 For each CQT Hopf $\V$-face algebra ${\HH}$, we have\rom{:}
 \begin{equation}
 \label{S(U)}					
  S^*({\cal U}_1^{\pm 1}) = {\cal U}_2^{\mp 1}, 
 \quad
  S^*({\cal U}_2^{\pm 1}) = {\cal U}_1^{\mp 1},  
 \end{equation}
 \begin{equation}
 \label{UXU-}
  {\cal U}_{\nu} X {\cal U}_{\nu}^{-1} =(S^2)^*(X)
  \quad (X \in {\HH}^*,\, \nu = 1, 2).
 \end{equation}
 In particular, $S$ is bijective and 
 ${\cal U}_1{\cal U}_2^{-1}$ is a central element of ${\HH}^*$.
 %
\end{prop}
\begin{pf}
 (cf. Drinfeld \cite{Drinfeld}).
 The relation \eqref{S(U)} follows from \eqref{Qp=Rm}, \eqref{Sid(R)}
 and \eqref{U-}.						
 Substituting $\sumc c_{(2)}$ $\otimes$ $S(c_{(1)})$ into
 $\Rp m^* (X) = (m^{\op})^* (X) \Rp$, 
 we obtain 
 \begin{align*}
  \sumc \cal{U}_1 (c_{(2)}) c_{(3)} S(c_{(1)})
  = & \sumc S(c_{(2)}) c_{(3)} \Qm (c_{(4)},c_{(1)}) \\
  = & \sumc \sumk e_k \Qm (c_{(2)} {\emaruk},c_{(1)}) \\
  = & \sumk \cal{U}_1 ({\emaruk}c) e_k,
 \end{align*}
 where the second equality follows from \eqref{e(eae)a}
 and the third equality follows from \eqref{Rmeeee} and \eqref{D(eeaee)}.
 Using this relation, we compute
 \begin{align*}
  \sumc S^2 (c_{(1)}) \cal{U}_1 (c_{(2)})  
  = & \sumc \sumk \cal{U}_1 ({\emaruk}c_{(2)}) e_k S^2(c_{(1)}) \\
  = & \sumc \cal{U}_1 (c_{(3)}) c_{(4)} S(S(c_{(1)}) c_{(2)}) \\
  = & \sumc \sumk \cal{U}_1 (c_{(1)}{\emaruk}) c_{(2)}{\emaruk} \\
  = & \sumc \cal{U}_1 (c_{(1)}) c_{(2)},  
  \end{align*}
 where the first equality follows from \eqref{D(a)} and \eqref{S(ee)}
 and the last equality follows from \eqref{U(ee)} and \eqref{D(a)}.
 Substituting this into $X \in {\HH}^*$, we get
 $((S^2)^*(X) \cal{U}_1)(c) = (\cal{U}_1 X)(c)$, 
 which proves \eqref{UXU-} for $\nu = 1$.
\end{pf}

%
%
%

\section{Group-like functionals}
   
Let $g$ be an element of a 
$\V$-face algebra $\HH$.
We say that $g$ is {\it group-like}
if
\begin{equation}
\label{D(g)}		
 \Delta (g) = \sumk g e_k \otimes g \emak,
\end{equation}
\begin{equation}
\label{gee}
 g \emai e_j = \emai e_j g, 
\quad
 \cu (g \emai e_j) = \delta_{ij}
\end{equation}
for each $i,j \in \V$.
We say that a linear functional 
${\cal G}$ on $\HH$ is 
{\it group-like}
if it is group-like as an element of the dual
face algebra $\HH^{\circ}$.
Explicitly, ${\cal G}$ is group-like
if and only if it satisfies
\begin{equation}
\label{Gab}
 {\cal G}(ab) = \sumk
 {\cal G} (ae_k) {\cal G}(\emak b),
\end{equation}			
\begin{equation}
\label{Geae}
 {\cal G}({\emarui}a{\emaruj}) = {\cal G}(e_i a e_j),
\end{equation}
\begin{equation}
\label{Gee}
 {\cal G}({\emarui}e_j) = \delta _{ij} 
\end{equation}
for each $a,b \in \HH$ and $i,j \in {\V}$.
We say that $\cal{G}$ is {\it invertible} if it is 
invertible as an element of the dual algebra $\HH^*$.
We denote by $\GLF (\HH)$
the set of all group-like functionals 
of $\HH$, and by $\GLF (\HH)^{\times}$
the set of all invertible group-like functionals.
Note that 
\begin{equation}
 \GLF (\HH) = \mathrm{Hom}_{\, \Bbb{K}\! - \! Alg} (\HH, \Bbb{K})
\end{equation}
if $\HH$ is a bialgebra.
\begin{lem}
\label{GLFisfunc}
\rom{(1)}
 The correspondence $\HH \mapsto \GLF (\HH)$ defines 
 a contravariant functor from the category of
 $\V$-face algebras to the category of semigroups. \\
\rom{(2)}				
 Let $\HH$ be a $\V$-face algebra and $\frak{I}$ its biideal.
 Then the projection $p\!: \HH \to \KK = \HH / \frak{I}$
 gives 
 \begin{equation}
  p^*\!: \GLF( \KK ) \cong 
  \{ \cal{G} \in \GLF( \HH )\, |\, \cal{G} (\frak{I}) = 0 \}.
 \end{equation} 
\rom{(3)}
 If $\HH$ has an antipode, then
 $\GLF (\HH) = \GLF (\HH)^{\times}$ 
 and 
 \begin{equation}
 \label{S(G)}
 S^*({\cal G}) = {\cal G}^{-1}. 
 \end{equation} 
 for each ${\cal G} \in \GLF (\HH)$. 
\end{lem}
\begin{pf}
 The proof of Part (1) is straightforward.
 Taking the dual of $0 \to \frak{I} \to \HH \to \KK \to 0$,  
 we obtain
 \begin{equation}
  p^*\!: \KK^* \cong 
  \{ X \in \HH^*\, |\, X (\frak{I}) = 0 \}.
 \end{equation} 
 It is straightforward to verify that 
 $\cal{M} \in \GLF (\KK)$ if and only if 
 $p^* (\cal{M}) \in \GLF (\HH)$ for each 
 $\cal{M} \in \KK^*$.
 This proves Part (2).
 See \cite{cpt} Proposition 7.1 for a proof of Part (3).
\end{pf}
\begin{lem}
\label{cqtglf}	
 Let $\HH$ be a CQT $\V$-face algebra. \\
 \rom{(1)} 
 For each  group-like functional $\cal{G}$ on $\HH$,
 we have 
 \begin{equation}
 \label{GGRpm}
  (\cal{G} \otimes \cal{G}) \Rpm = 
  \Rpm (\cal{G} \otimes \cal{G}).
 \end{equation} 
 Hence, for each $\HH$-comodules $U$ and $V$, we have  
 \begin{equation}
 \label{GGcMN}					
  (\pi_{V} (\cal{G}) \bar{\otimes} \pi_{U} (\cal{G})) c_{UV} = 
  c_{UV} (\pi_{U} (\cal{G}) \bar{\otimes} \pi_{V} (\cal{G})).
 \end{equation}   
 \rom{(2)} If $\HH$ is closable, then 
 \begin{equation}
 \label{U1U2isglf}
 \cal{U}_1  \cal{U}_2 \in
 \GLF (\HH).
 \end{equation}
\end{lem}
\begin{pf}
 Since $\Rp = (m^{\op})^* (1) \Rp$,
 we have			
 $(\cal{G} \otimes \cal{G}) \Rp$ $=$
 $(m^{\op})^* (\cal{G}) \Rp$.
 Hence the first assertion of Part (1) follows from
 \eqref{RmXR}.
 The second assertion follows from the first assertion.
 Part (2) follows from \eqref{U(ee)} and \eqref{m(U)}-\eqref{m(Um)}.
\end{pf}
Let $\HH$ be a $\V$-face algebra and $\cal{G}$ its group-like functional.
We define $\mathrm{coad}(\cal{G})\!: \HH \to \HH$ by
\begin{equation}
 \label{coaddef}
 \mathrm{coad}(\cal{G}) (a) =
 \suma \cal{G}^{-1} (a_{(1)}) a_{(2)} \cal{G} (a_{(3)})
 \quad (a \in \HH).
\end{equation}
Using \eqref{Gab}-\eqref{Gee} and \eqref{eae*a},
we  see that $\mathrm{coad}(\cal{G})$ is 
an automorphism of $\HH$.
\begin{prop}
\label{GLFAwiso}
For $\HH =$ $\HG$ or $\Aw$, the map 
$\HH^* \to \End (\K \G^1);$
$X \mapsto \pi_{\K \G^1} (X)$
gives the following semigroup isomorphisms\rom{:}
\begin{equation}
 \label{GLFHG}
 \GLF (\HG) \cong 
 \End_{\pi (\frak{E})} (\K \G^1), 
\end{equation}
\begin{equation}
 \label{GLFAw}			
 \GLF (\Aw) \cong \bigl\{ G \in \End_{\pi (\frak{E})} (\K \G^1) \bigm|  
 (G \bar{\otimes} G) w =
 w (G \bar{\otimes} G)
 \bigr\},
\end{equation}
where $\frak{E} = \frak{E}_{\HH^{\circ}}$ is as in Sect. 2.
\end{prop}
\begin{pf}
For each element $G$ of the right-hand side of \eqref{GLFHG}, 
we define a linear functional $\cal{G}$ $=$ $\cal{G}^{\HH}_G$ 
on $\HG$ by setting
\begin{equation}
\label{Gepq}
 \cal{G} \left( e \binom{i}{j} \right) = \delta_{ij}, 
\quad
 \cal{G} \left( e \binom{\p}{\q} \right) = 
 G^{\p_1}_{\q_1} 
 \cdots G^{\p_m}_{\q_m} 
\end{equation}
for each paths $\p = (\p_1, \ldots, \p_m)$ and
$\q = (\q_1, \ldots, \q_m)$ of length $m > 0$ and 
$i, j \in \V$.
It is straightforward to verify that $\cal{G}$ is 
a group-like functional of $\HG$.
Hence $\pi_{\K \G^1}$ gives a surjection
$\GLF (\HG) \to \End_{\pi (\frak{E})} (\K \G^1)$. 
Conversely, for $\cal{G} \in \GLF (\HG)$, 
set $G = \pi_{\K \G^1} (\cal{G})$.
Then by \eqref{Gab}-\eqref{Gee}, we have
\eqref{Gepq}.
Thus we get the isomorphism \eqref{GLFHG}.
Next we show \eqref{GLFAw}.
By \eqref{GGcMN}, $\pi_{\K \G^1}$ defines a 
well-defined map from $\GLF (\Aw)$ to the right-hand side of
\eqref{GLFAw}.
Hence it suffices to construct the inverse of this map.
Let $G$ be an element of the right-hand side of 
\eqref{GLFAw} and let $\cal{G}^{\HH}_G \in \GLF (\HG)$ be as above.			
By \eqref{GGcMN}, we have
\begin{equation}
\label{GrelAw}
 \cal{G}^{\HH}_G \left( 
 \sum_{(\ccc, \ddd) \in {\G}^2}
 w \!\! \left[ \aaa \frac[0pt]{\ccc}{\bbb} \ddd \right] 
 e{\ccc \cdot \ddd  \choose \p \cdot \q} -
 \sum_{(\r, \s) \in {\G}^2}
 w \!\! \left[ \r \frac[0pt]{\p}{\s} \q \right]
 e{\aaa \cdot \bbb  \choose \r \cdot \s}
 \right) = 0
\end{equation}
for each $(\p, \q),\, (\aaa, \bbb) \in {\G}^2$.
By \eqref{Gab}, this shows that 
$\cal{G}^{\HH}_G$ vanishes on
the biideal $\mathrm{Ker} (\HG \to \Aw)$ and that
it induces an element of $\GLF (\Aw)$. 
This completes the proof of \eqref{GLFAw}.
\end{pf}
\begin{prop}
\label{GLFHcHHiso}
 For each CCQT $\V$-face algebra, 
 the canonical map 
 $\iota\!: \HH \to \HcH$ induces the isomorphism
 \begin{equation}
  \iota^*\!: \GLF (\Hc (\HH)) \cong 
  \GLF^{\times} (\HH), 
 \end{equation}
 whose inverse $\cal{G} \mapsto \cal{G}_{\Hc}$ is given by
 \begin{equation}
 \label{Gasb}
  \cal{G}_{\Hc} (a \sigma (b)) = 
  \sumk \cal{G}(a e_k) \cal{G}^{-1} (b e_k).
 \end{equation}
\end{prop}
\begin{pf}%
 By Lemma \ref{GLFisfunc} (1), it suffices to show that \eqref{Gasb}
 gives the inverse of the correspondence $\iota^*$.	
 It is easy to verify that there exists a linear functional 
 $\hat{\cal{G}} \in \hat{\HH}^*$ which sends 
 $a \tensE \sigma (b)$ to the right-hand side of \eqref{Gasb}
 and that $\hat{\cal{G}}$ satisfies \eqref{Geae} and \eqref{Gee}.
 Using \eqref{Gab} for $\cal{G}^{\pm 1}$, we obtain
 \begin{equation}
  \hat{\cal{G}} ((a \tensE 1) x
  (1 \tensE \sigma (d))) =
  \sum_{i,j \in \V} \cal{G} (a e_i) 
  \hat{\cal{G}} (\emai x e_j)
  \cal{G}^{-1} (d e_j) 
  \quad (a, d \in \HH, x \in \hat{\HH}).
 \end{equation}
 By replacing $x$ with 
 $(1 \tensE \sigma (b))(c \tensE \sigma(1))$, we obtain 
 \begin{multline}
 \label{Gasbcsd}
  \hat{\cal{G}} ( (a \tensE \sigma (b))(c \tensE \sigma(d)) ) 
  = \sum_{i,j \in \V} \cal{G} (a e_i) 
 \hat{\cal{G}} ( (1 \tensE \sigma (b e_i))(c e_j \tensE \sigma(1)) )  
 \cal{G}^{-1} (d e_j).
 \end{multline}
 On the other hand, using \eqref{Geae},
 \eqref{eae*a} and \eqref{Rmeeee},
 we obtain
 \begin{multline}
  \hat{\cal{G}} ( (1 \tensE \sigma (b))(c \tensE \sigma(1))) = 
  \sumk \sumbc \Rm (b_{(1)} e_k , c_{(3)} \emak)
  \Qp (b_{(3)}, c_{(1)})
  \cal{G}^{-1} (b_{(2)}) \cal{G} (c_{(2)})  \\
  = \sumbc 
  \langle (1 \otimes \cal{G}) \Rm (\cal{G}^{-1} \otimes 1),\,
  b_{(1)} \otimes c_{(2)} \rangle
  \Qp (b_{(2)}, c_{(1)}).
 \end{multline}
 By \eqref{GGRpm},
 the right-hand side of the above equality is
 \begin{multline}
  \sumbc 			
  \cal{G}^{-1} (b_{(1)})
  \Rm ( b_{(2)}, c_{(2)}) \Qp (b_{(3)}, c_{(1)}) 
  \cal{G} (c_{(3)})
  = \sumk \cal{G}^{-1} (e_k b) \cal{G} (e_k c).
 \end{multline}
 Hence the right-hand side of \eqref{Gasbcsd} is 
 \begin{equation}
  \sum_{i,j,k \in \V}	
  \cal{G} (a e_i) \cal{G}^{-1} (e_k a e_i)
  \cal{G} (e_k c e_j) \cal{G}^{-1} (d e_j)
  = \sumk			
  \hat{\cal{G}} ((a \tensE \sigma (b)) e_k)
  \hat{\cal{G}} (e_k (c \tensE \sigma (d))).
 \end{equation}
 Thus $\hat{\cal{G}}$ is a group-like functional of $\hat{\HH}$.
 Using \eqref{Gab} for $\hat{\cal{G}}$, we compute 
 \begin{multline}
  \hat{\cal{G}} \left( \emai e_j \left( 
  \suma (1 \tensE \sigma (a_{(1)})) (a_{(2)} \tensE 1) 
 \right) \emak e_l \right) \\ 
  =
  \sum_{m \in \V} \suma \cal{G}^{-1} (\ema_m 
  a_{(1)} \emaj e_i ) 
  \cal{G} (\ema_m a_{(2)} \emak e_l) 
  = 
  \suma \delta_{ij} \delta_{kl}
  \cal{G}^{-1} (a_{(1)} \emaj) 
  \cal{G} (a_{(2)} e_l) \\
  = 
  \delta_{ij} \delta_{kl} \delta_{jl}
  \cu (a e_l)
  =
  \hat{\cal{G}} \left( \emai e_j \left( 
  \sum_{m \in \V} \cu (a e_m) e_m 
  \right) \emak e_l \right)  
  \qquad\qquad\qquad
 \end{multline}
 for each $i,j,k,l \in \V$ and $a \in \HH$,  
 where the second equality follows from 
 \eqref{Geae} and \eqref{D(a)} and
 the third equality follows from 
 \eqref{D(eeaee)}.
 By repeating similar calculation, we see that 
 $\hat{\cal{G}}$ induces a group-like functional 
 $\cal{G}_{\Hc}$ on $\HcH$.
 Now it is straightforward to verify that 
 $\cal{G} \mapsto \cal{G}_{\Hc}$ gives the inverse of $\iota^*$.
\end{pf}
Combining Proposition \ref{GLFAwiso} and 
Proposition \ref{GLFHcHHiso},
we obtain the group isomorphism
\begin{equation}
\label{Phidef}
 \Phi\!:
 \bigl\{ G \in \Aut_{\pi (\frak{E})} (\K \G^1) \bigm|  
 (G \bar{\otimes} G) w =
 w (G \bar{\otimes} G)\, \bigr\}  
 \cong
 \GLF (\Hc (\Aw))
\end{equation}
for each star-triangular face model $(\G, w)$.

\section{A classification theory of ribbon functionals}

\begin{lem}
 For a coribbon Hopf $\V$-face algebra ${\HH}$,
 we have
 \begin{equation}
 \label{V(eae)}			
  \cal{V}^{\pm 1} (\emai a \emaj) =   
  \cal{V}^{\pm 1} (e_i a e_j),
 \quad
  \cal{V}^{\pm 1} (\emai e_j) = \delta_{ij},
 \end{equation}
 \begin{equation}
 \label{V2}
  \cal{V}^2 = \cal{U}_1 \cal{U}_2^{-1},
 \end{equation}
 \begin{equation}
 \label{m(V-)}			
  m^* (\cal{V}^{-1}) = 
 (\cal{V} \otimes \cal{V})^{-1}
 \Rp_{21} \Rp.
 \end{equation}
\end{lem}
\begin{pf}
 The first equality of \eqref{V(eae)} follows from the fact that
 $\cal{V}$ commutes with the face idempotents of $\HH^{\circ}$.
 Using \eqref{m(V)} and \eqref{Rmeea}, we obtain
 \begin{align}
  \cal{V} (a) & = 
  \sum_{j,k \in \V} \sum_{(a)} \Rm (e_j, a_{(1)}) \Rm (a_{(2)}, \emaj e_k)
  \cal{V} (\emak) \cal{V} (a_{(3)})\\
  & = \sumk \cal{V} (\emak) \cal{V} (\emak a) \nonumber\\
  & = \langle 
  \sumk \cal{V} (\emak) \emak \cal{V}, 
  a \rangle, \nonumber
 \end{align}
 which implies
 $\sum_k \cal{V} (\emak) \emak = 1$.
 Since $\{ \emak \}$ is linearly independent by 
 the second equality of \eqref{D(ee)},
 this proves $\cal{V} (\emak) = 1$, 
 or the second equality of \eqref{V(eae)}.
 By a similar discussion to \cite{Kassel} page 351, 
 we obtain $\pi_M (\cal{V}^2) =$ 
 $\pi_M (\cal{U}_1 \cal{U}_2^{-1})$ for 
 every $\HH$-comodule $M$.
 Hence \eqref{V2} follows from the fundamental theorem
 of coalgebras (cf. \cite{Sweedler} page 46).
 Using the fact that $\Rm$ is the 
 $(m^* (1), (m^{\op})^* (1))$-generalized inverse of $\Rp$,
 we obtain
 \begin{equation}
 \label{RRRR}			  
  \Rm \Rm_{21} \Rp_{21} \Rp
  = m^* (1) = \Rp_{21} \Rp \Rm \Rm_{21}.
 \end{equation}
 Hence the right-hand side of \eqref{m(V-)} is 
 the inverse of $m^* (\cal{V})$ in the algebra 
 $m^* (1) (\HH^{\otimes 2})^*$ $m^* (1)$.
 This proves \eqref{m(V-)}.
\end{pf}
\begin{prop}
 Let ${\HH}$ be a CQT Hopf $\V$-face algebra and
 ${\cal V}$ an invertible element of $\HH^*$.  
 Then $(\HH,\cal{V})$ is a coribbon Hopf $\V$-face algebra 
 if and only if $\cal{M} = \cal{U}_1 \cal{V}^{-1}$ 
 is group-like and satisfies the following relations\rom:    
 \begin{equation}
 \label{MXM-}
 \qquad
  {\cal M} X {\cal M}^{-1} = 
  (S^2)^* (X)
  \quad (X \in \HH^*),
 \end{equation}
 \begin{equation}
 \label{M2=U1U2}
  \cal{M}^2 = \cal{U}_1 \cal{U}_2. \qquad
 \end{equation}
\end{prop}
\begin{pf}
 To begin with, we note that the equivalence of
 $\cal{V} \in Z(\HH^*)$ and \eqref{MXM-} follows from \eqref{UXU-},
 and that that of
 \eqref{m(V)} and \eqref{D(g)} for $g = \cal{M}$ follows from
 \eqref{m(U)}, \eqref{m(V-)} and \eqref{RRRR}.
 Suppose $\cal{V}$ is a ribbon functional.
 Then the relation \eqref{M2=U1U2} follows from \eqref{V2}
 and \eqref{U1U2},
 while the first (resp. second) relation of \eqref{gee} 
 for $g = \cal{M}$ follows from \eqref{UXU-}
 and \eqref{S(ee)} (resp. \eqref{U(ee)} 
 and the second relation of \eqref{V(eae)}).
 Conversely, if $\cal{M}$ satisfies the above conditions, 
 \eqref{S(V)} follows from \eqref{V2}, \eqref{U1U2} and \eqref{S(G)}.
\end{pf}
For a coribbon Hopf $\V$-face algebra $(\HH, \cal{V})$, 
we call $\cal{M} = \cal{U}_1 \cal{V}^{-1}$ the 
{\it modified ribbon functional} on $\HH$
corresponding to $\cal{V}$. 
For each CQT Hopf $\V$-face algebra $\HH$, we denote by     
$\mathrm{Rib}(\HH)$ the set of all ribbon functionals on $\HH$
and by $\mathrm{MRib}(\HH)$ the set of all 
modified ribbon functionals on $\HH$.
\begin{prop}
\label{MRibaltbraiding}
 Let $(\HH, \Rpm)$ be a CQT Hopf $\V$-face algebra.
 \rom{(1)}
 We have
  \begin{equation}
  \label{MRibR-21}
   \mathrm{MRib}((\HH, \Rmp_{21})) = 
   \mathrm{MRib}((\HH, \Rpm)).  
  \end{equation}
 \rom{(2)}
 Let $\HH_{\gamma}$ $(\gamma \in \Gamma)$
 and $\chi$ be as in Proposition \ref{altbraiding}.
 Then we have
 \begin{equation}
  \label{MRibRchi}
   \mathrm{MRib}((\HH, \Rpm_{\chi})) = 
   \mathrm{MRib}((\HH, \Rpm)).  
  \end{equation}
\end{prop}
\begin{pf}
 Let $\cal{U}_i$, $\cal{U}_i^{\prime}$  
 and $\cal{U}_{i, \chi}$ $(i = 1,2)$ be 
 the Drinfeld functionals of $(\HH, \Rpm)$, 
 $(\HH, \Rmp_{21})$ and
 $(\HH, \Rpm_{\chi})$ respectively.
 Then we have $\cal{U}_1^{\prime}$ $=$ $\cal{U}_2$,
 $\cal{U}_2^{\prime}$ $=$ $\cal{U}_1$ and 
 \begin{equation}	
  \cal{U}_{1, \chi} (a) = \chi (\gamma, \gamma)^{-1} \cal{U}_1 (a), 
 \quad
  \cal{U}_{2, \chi} (a) = \chi (\gamma, \gamma) \cal{U}_2 (a)
  \quad (a \in \HH_{\gamma}).
 \end{equation}
 Hence the assertions follows from the definition of 
 the modified ribbon functional and \eqref{U1U2}.
\end{pf}
\begin{thm}
\label{clasRibHcAw}				
 For each closable star-triangular face model $(\G, w)$,
 the map $\pi_{\K \G^1}$ gives the following bijection\rom{:}
 \begin{equation}
 \label{Ribiso}
  \mathrm{Rib}(\Hc (\Aw)) \cong 
  \{ V \in \mathrm{Aut}_{\Hc (\Aw)} (\K \G^1)\, |\, 
  V^2 = \pi_{\K \G^1} (\cal{U}_1 \cal{U}_2^{-1}) \}.
 \end{equation}
 Equivalently, $\pi_{\K \G^1}$ gives
 \begin{equation}
 \label{MRibiso}
  \mathrm{MRib}(\Hc (\Aw)) \cong 
  \{ M\, |\, M \pi_{\K \G^1} (\cal{U}_1)^{-1}
  \in \mathrm{Aut}_{\Hc (\Aw)} (\K \G^1),\, 
  M^2 = \pi_{\K \G^1} (\cal{U}_1 \cal{U}_2) \}. 
 \end{equation}
\end{thm}
\begin{pf}				
 Let $M$ be an element of the right-hand side of \eqref{MRibiso}.
 By \eqref{UURp}, 
 $\pi_{\K \G^1} (\cal{U}_{\nu}) \bar{\otimes}$
 $\pi_{\K \G^1} (\cal{U}_{\nu})$ 
 commutes with $w$ for each $\nu = 1, 2$.		
 Hence $M$ belongs to 
 the right-hand side of \eqref{GLFAw}.		
 Set $\cal{M} = \Phi (M)$, where $\Phi$ is as in \eqref{Phidef}.
 Since $\pi_{\K \G^1} (\cal{M}^2)$ $=$ 
 $\pi_{\K \G^1} (\cal{U}_1 \cal{U}_2)$, we have 
 $\cal{M}^2$ $=$ $\cal{U}_1 \cal{U}_2$ by 
 Lemma \ref{cqtglf} (3).
 By \eqref{UXU-}, we have 
 $\mathrm{coad}(\cal{M}) (e \binom{\p}{\q})$
 $=$ $S^{-2} (e \binom{\p}{\q})$,
 for each $\p, \q \in \G^1$. 
 Since $\mathrm{coad}(\cal{M})$ is an automorphism and 
 $e \binom{\p}{\q}$, $S(e \binom{\p}{\q})$
 $(\p, \q \in \G^1)$ and 
 $e \binom{i}{j}$ $(i,j \in \V)$ generate $\Hc (\Aw)$, 
 this shows that $\mathrm{coad}(\cal{M}) = S^{-2}$. 
 Thus $\cal{M}$ is a modified ribbon functional of
 $\Hc (\Aw)$. Conversely, it is clear that 
 $\pi_{\K \G^1}$ maps
 the left-hand side of \eqref{MRibiso}
 into the right-hand side of \eqref{MRibiso}.
 Thus we get the theorem.
\end{pf}		
\begin{thm}[\cite{Reshetikhin}]		
\label{existrib}
 For each closable star-triangular face model $(\G, w)$
 over an algebraically closed field $\K$ of $\mathrm{ch} \K \ne 2$, 
 $\Hc (\Aw)$ has a ribbon functional.
\end{thm}
\begin{pf}				
 It suffices to construct a linear operator $V$
 which belongs to the right-hand side of \eqref{Ribiso}.
 Let $A$ be the operator $\pi_{\K \G^1} (\cal{U}_1 \cal{U}_2^{-1})$
 and $A = S + N$ its Jordan decomposition, that is,
 $S$ is a diagonalizable operator and $N$ is a nilpotent operator
 such that $SN = NS$.
 Let $\lambda_i$ $(1 \leq i \leq k)$ be (mutually distinct)
 eigenvalues of $S$ and $P_i$ the projection corresponding to 
 $\lambda_i$.
 It is known that $P_i = f_i (A)$ and   
 $N = g(A)$ for some polynomials $f_i, g \in \K [X]$.
 Let $ \sqrt{\lambda_i}$ be a square root of $\lambda_i$
 and define a operator $V$ by
 $V = \sum_i \sqrt{\lambda_i} P_i h( S^{-1} N)$, 
 where $h \in \K [X]$ is defined by 
 \begin{equation} 
  h(X) = 1 + \sum_{n = 0}^{\sharp \G^1} (-1)^n 2^{-2n - 1}
  \frac{1}{n + 1} \binom{2n}{n} X^{n + 1}.
 \end{equation} 
 Then, we have  $V^2 = A$.
 Since $\cal{U}_1 \cal{U}_2^{-1}$ is a central element of 
 $\Hc (\Aw)^*$ and $V$ is a polynomial of $A$, we have
 $V \in \mathrm{Aut}_{\Hc (\Aw)} (\K \G^1)$. 
 By the theorem above, this proves the existence of a
 ribbon functional on $\Hc (\Aw)$.
\end{pf}				
Let $(\G, w)$ be a closable star-triangular face model.
We say that $(\G, w)$ is ({\it absolutely}) 
{\it irreducible} if 
$\K \G^1$ is (absolutely) irreducible as an $\HcH$-comodule.
As an immediate consequence of the Theorem \ref{clasRibHcAw} and
Schur's Lemma, we have the following.
\begin{thm}
\label{cardRib}		
 Let $(\G, w)$ be an irreducible
 closable star-triangular face model
 over an algebraically closed field.
 Then we have
 $\sharp \mathrm{Rib}(\Hc (\Aw)) = 2$ if $\mathrm{ch} \K \ne 2$
 and 
 $\sharp \mathrm{Rib}(\Hc (\Aw)) = 1$ if $\mathrm{ch} \K = 2$.
\end{thm}
\begin{thm}				
\label{Mcrit}
 Let $(\G,w)$ be an absolutely irreducible closable
 star-triangular face model.
 Suppose $M \in GL (\K \G^1)$ satisfies 
 $\sum_{\r \s} M^{\p}_{\,\r} e \binom{\r}{\s} 
 (M^{-1})^{\s}_{\q}$ $=$ $S^2 (e \binom{\p}{\q})$ and
 $\mathrm{Tr} (M) =  \mathrm{Tr} (M^{-1}) \ne 0$.
 Then we have 
 \begin{equation} 
  \mathrm{MRib} (\Hc (\Aw)) = \{ \Phi (\pm M) \}.
 \end{equation} 
\end{thm}
\begin{pf}
 By Schur's lemma, we have 
 $\pi (\cal{U}_{\nu})$ $=$ $c_{\nu} M$
 for some nonzero constant $c_{\nu}$ $(\nu = 1,2)$.
 Since 
 $\mathrm{Tr} \pi (\cal{U}_1)$ $=$
 $\mathrm{Tr} \pi (\cal{U}_2^{-1})$
 by \eqref{Udef} and \eqref{U-}, 
 we obtain  
 $c_1 \mathrm{Tr} (M)$ $=$
 $c_2^{-1} \mathrm{Tr} (M^{-1})$.
 Therefore $M$ belongs to the right-hand side of \eqref{MRibiso}.
\end{pf}
\begin{prop}
\label{clasquotient}
 Let $(\G, w)$ be a closable star-triangular face model 
 and let $\KK = \Hc (\Aw) / \frak{I}$ be a quotient CQT 
 Hopf $\V$-face algebra of $\HH:= \Hc (\Aw)$ such that 
 $\K \G^1$ is absolutely irreducible as a $\KK$-comodule.
 Then the projection $p\!: \HH \to \KK$ gives the isomorphism
 \begin{equation} 
 \label{RibH}
  p^*\!: \mathrm{MRib} (\KK) \cong
  \{ \cal{M} \in \mathrm{MRib} (\Hc (\Aw))\, |\, \cal{M}(\frak{I}) = 0 \}. 
 \end{equation} 
\end{prop}
\begin{pf}
 We prove the assertion by using 
 Lemma \ref{GLFisfunc} (2). 
 Let $\cal{M}$ be a group-like functional on $\KK$.
 It suffices to verify that 
 $\mathrm{coad}(\cal{M})$ $=$ $S^{-2}$ if and only if 
 $\mathrm{coad}(p^* (\cal{M}))$ $=$ $S^{-2}$.
 Since $\mathrm{coad}(\cal{M})(p (a))$ $=$ 
 $p (\mathrm{coad}(p^* (\cal{M}))(a))$
 for each $a \in \HH$, the ``if''-part is 
 obvious.				
 Suppose  $\mathrm{coad}(\cal{M})$ $=$ $S^{-2}$ and 
 set $M: = \pi_{\K \G^1} (\cal{M})$.
 Since $\pi_{\K \G^1}^{\HH} (p^* (\cal{M})) = M$, 
 we have $\Phi (M) = p^* (\cal{M})$. 
 On the other hand, using  \eqref{f(U)=U} and Schur's Lemma, 
 we see that $M \pi_{\K \G^1}^{\HH} (\cal{U}_1)^{-1}$
 is a scalar multiple of the identity operator.
 Hence $M$ belongs to the right-hand side of \eqref{MRibiso}. 
 By Theorem \ref{clasRibHcAw}, this proves the proposition.
\end{pf}
Let $\HH$ be a CQT Hopf $\V$-face algebra.
We say that $\HH$ is {\em monogenerated}
if there exists an absolutely irreducible 
$\HH$-comodule $U$ such that $\HH$ is generated by
$\emai e_j$ $(i,j \in \V)$, the image $C$ of the 
corepresentation $\End(U)^* \to \HH$ and $S(C)$,
as an algebra.
\begin{lem} 
 Let $\HH$ be a CQT Hopf $\V$-face algebra over $\K$ and 
 $\Bbb{F}$ a field extension of $\K$.
 Then $\HH \otimes \Bbb{F}$ naturally
 becomes a CQT Hopf $\V$-face algebra over $\Bbb{F}$ and
 there exists an injection 
 $\mathrm{Rib} (\HH) \to \mathrm{Rib} (\HH \otimes \Bbb{F});$
 $\cal{V} \mapsto \cal{V}_{\Bbb{F}}$
 given by 	
 $\cal{V}_{\Bbb{F}} (a \otimes 1_{\Bbb{F}})$
 $=$ $\cal{V} (a)$. 
\end{lem}		
\begin{pf}
 This is straightforward.
\end{pf}
\begin{prop}		
\label{ribestimate}		
 For each monogenerated CQT Hopf $\V$-face algebra $\HH$,
 we have
 $\sharp \mathrm{Rib}(\HH) \leq 2$ if $\mathrm{ch} \K \ne 2$
 and
 $\sharp \mathrm{Rib}(\HH) \leq 1$ if $\mathrm{ch} \K = 2$
\end{prop}			
\begin{pf}
 Let $(\G, w_U)$ and
 $f\!: \Hc (\frak{A}(w_U)) \to \HH$ be as in Proposition \ref{univA}.
 Since $\HH$ is monogenerated, $f$ is surjective
 for a suitable absolutely irreducible comodule $U$.  
 Now the assertion is an immediate consequence of 
 \eqref{RibH}, Theorem \ref{cardRib}
 and the lemma above.
\end{pf}
\noindent
{\it Remark.}
 (1)
 To construct a link invariant via a lattice model $(\G, w)$, 
 it is usual to assume that $(\G, w)$ is ``enhanced'' in the sense
 of \cite{TuraevYB} (cf. \cite{ADW}, \cite{Jones}, \cite{TuraevYB}). 
 Theorem \ref{existrib} says that the assumption is superfluous 
 provided that $(\G, w)$ is closable.
 For vertex models, this was first proved by 
 Reshetikhin \cite{Reshetikhin}. \\
 (2) 
 Combining Theorem \ref{existrib} with the categorical framework
 of the link invariant \cite{Turaev3mfd}, we obtain an invariant of
 framed links colored by comodules of $\Hc (\Aw)$, 
 for each closable star-triangular face model $(\G, w)$.
 Choosing the $\Hc (\Aw)$-comodule $\K \G^1$ as a color, 
 we obtain an invariant $I_w (L)$ of framed links $L$
 which agrees with the known one.   
 However, if $(\G, w)$ is constructed 
 from a (four-weight) spin model $(W_i)$ (\cite{Jones}, \cite{Bannai^2}), 
 $I_w (L)$ does not agree with the known invariant $Z_{(W_i)} (L)$.
 In fact we have $I_w (L) = Z_{(W_i)} (L) Z_{(W_i)}^* (L)$, where 
 $Z_{(W_i)}^* (L)$ is the ``dual invariant'' of $Z_{(W_i)} (L)$.

\section{Quantized classical groups}

%
%
%
Let $X_l$ be one of the Dynkin  diagram of type $A_l$,
$B_l$, $C_l$ or $D_l$, where $l \geq 1$ if $X = A$ and 
$l \geq 2$ if $X = B, C$ or $D$.
We define integers $N$ and
$\nu$ by
\begin{equation}
N = 
 \begin{cases}
  l + 1  & (X = A) \\
  2l + 1 & (X = B) \\
  2l     & (X = C, D), \\
 \end{cases}
\quad
 \nu = 
 \begin{cases}
  0  & (X = A) \\
  - 1 & (X = B, D) \\
  1   & (X = C).
 \end{cases}
\end{equation}
For $X = B, C, D$ and $1 \leq i \leq N$, 
we set $i^{\prime} = N + 1 - i$ and  
$\bar{i} = i - \sigma_i \nu / 2$, where
\begin{equation}
 \sigma_i = 
 \begin{cases}
  1   & (1 \leq i < (N+1)/2) \\
  0   & (i = (N+1)/2) \\
  - 1 &((N+1)/2 < i \leq N),
 \end{cases}
\quad
 \epsilon_i =
  \begin{cases}
  1     & (1 \leq i \leq (N+1)/2) \\
  - \nu & ((N+1)/2 \leq i \leq N).
 \end{cases}
\end{equation}
Also we set $\sigma_i \equiv 1$ for $X = A$.
Let $\Check{R} = \Check{R}_q (X_l)$ 
be Jimbo's solution of 
the Yang-Baxter equation of type $X_{l}$:
\begin{equation}
\label{RAdef}
 \Check{R}_q (A_l) =
 q^{-1} \sum_{r = 1}^{N}
 E_{rr} \otimes E_{rr} +
 \sum_{r \ne s}
 E_{rs} \otimes E_{sr} -
 (q - q^{-1}) \sum_{r > s}
 E_{rr} \otimes E_{ss}, 
\end{equation}
\begin{multline}
\label{RBCDdef}
 \Check{R}_q (X_l) =
 \sum_{r;\, r \not= r'}
 (q^{-1} E_{rr} \otimes E_{rr} + 
 q E_{rr'} \otimes E_{r'r}) + 
 \sum_{r;\, r = r'}
 E_{rr} \otimes E_{rr} + \\
 \sum_{r,s;\, r \ne s,s'}
 E_{rs} \otimes E_{sr} +
 (q - q^{-1}) \sum_{r > s}
 (- E_{rr} \otimes E_{ss} 
 + \epsilon_r \epsilon_s q^{\overline{r} - \overline{s}}
 E_{rs'} \otimes E_{r's})
 \quad (X = B, C, D), 
\end{multline}
where for $X =$ $A, C, D$ (resp. $X= B)$, $q$ 
(resp. $q^{1/2}$) denotes
a non-zero number such that $q^2 \not= 1$,
and
$E_{rs} \in \mathrm{Mat} (N, \Bbb{K})$
denote the matrix units.
For $X =$ $B, C, D$,  we also set $\lambda = - \nu q^{-N - \nu}$.

For $1 \leq i, j \leq N$, we denote by
$t_{ij}$ the element 
$e {i \choose j}$ of $\frak{A} (\Check{R})$, or its image by an arbitrary
bialgebra map.
For each $\eta \in \K^{\times}$, we denote by $\Rp_{\eta}$
the canonical braiding of the FRT bialgebra 
$\frak{A} (\eta \Check{R})$ or its Hopf closure
$\mathrm{Hc} (\frak{A} (\eta \Check{R}))$.
Since $\frak{A} (\Check{R})$ 
(resp. $\mathrm{Hc} (\frak{A} (\Check{R}))$)
is isomorphic to
$\frak{A} (\eta \Check{R})$ 
(resp. $\mathrm{Hc} (\frak{A} (\eta \Check{R}))$)
as a bialgebra, 
we regard $\{ \Rp_{\eta} \}$ as a one-parameter family of 
braidings of $\frak{A} (\Check{R})$
(resp. $\mathrm{Hc} (\frak{A} (\Check{R}))$).

\begin{thm}[Takeuchi \cite{Takeuchicocycle}]
\label{clasbrFRT}	
 Any braidings of $\frak{A} (\Check{R}_q (X_{l}))$
 are either of the form $\Rp_{\eta}$ 
 or of the form $(\Rm_{\eta})_{21}$, 
 where $\eta \in \K^{\times}$.
\end{thm}
\begin{pf}
 For $X = A$, this theorem has been proved by  
 M. Takeuchi \cite{Takeuchicocycle}.
 Here we give a proof for $X = B, C, D$
 by imitating his arguments. 
 It is well known that the operators
 $g\!:= \Check{R}$ and
 $e\!:= (g - g^{-1}) / \mu + 1$
 give a representation of the Birman-Murakami-Wenzl algebra on 
 $(\K^N)^{\otimes 3}$ (cf. \cite{BirmanWenzl}, \cite{Murakami}), 
 where $\mu = q - q ^{-1}$.
 That is, we have the following formulas:
 \begin{equation}
 \label{minpolyR} 
 (g_i - \lambda^{-1})(g_i + q )(g_i - q^{-1}) = 0
 \quad (i = 1, 2),
 \end{equation}
 \begin{equation}
   g_1 g_2 g_1 = g_2 g_1 g_2,
 \quad
  e_1 g_2 e_1 = \lambda e_1,  
 \quad
  e_2 g_1 e_2 = \lambda e_2,  
 \end{equation}
 where, as usual, we set
 $f_1$ $=$ $f \otimes \mathrm{id}_{\K^N}$ and
 $f_2 = \mathrm{id}_{\K^N} \otimes f$ 
 for $f \in \End_{\K} ((\K^N)^{\otimes 2})$.
 As consequences of these relations, we also obtain 
 the following formulas:
 \begin{equation}
 \label{gi2} 
  g_i^2 = - \mu g_i + \lambda^{-1} \mu e_i + 1, 
 \quad 
  e_i^2 = \zeta e_i, 
 \end{equation}						
 \begin{equation}
 \label{eigi} 
  e_i g_i = g_i e_i = \lambda^{-1} e_i,  
 \end{equation}
 \begin{equation}
  e_i e_j e_i = e_i,  
 \quad 
  e_i g_j e_i = \lambda e_i,  
 \end{equation}
 \begin{equation}
  e_i e_j g_i = e_i g_j - \mu e_i e_j + \mu e_i,
 \quad 
  g_i e_j e_i = g_j e_i - \mu e_j e_i + \mu e_i,
 \end{equation}
 \begin{equation}
  e_i g_j g_i = e_i e_j,
 \quad 
  g_i g_j e_i = e_j e_i,
 \end{equation}
 \begin{equation}
  g_i e_j g_i - 
  g_j e_i g_j = 
  \mu (e_i g_j + g_j e_i - e_j g_i - g_i e_j) +
  \mu^2 (e_i - e_j)
 \end{equation}
 for $(i,j)$ $=$ $(1,2)$, $(2,1)$,
 where $\zeta = $ $- (\lambda - \lambda^{-1}) \mu^{-1} + 1$.
 By \eqref{gi2} and \eqref{eigi}, we see that 
 $\{ g, e, 1 \}$ is a linear basis of the algebra
 $\langle g \rangle$.

 Let $\cal{B}$ be a braiding of $\frak{A} (\Check{R}))$
 and $\Check{B} \in \End ((\K^N)^{\otimes 2})$ 
 the corresponding solution of the Yang-Baxter equation.
 Since $\frak{A}_2 (\Check{R})^*$ is the commutant of the algebra
 $\langle g \rangle$
 in $\End_{\K} ((\K^N)^{\otimes 2})$, 
 $\Check{B}$ belongs to the double commutant of 
 $\langle g \rangle$.
 By \cite{Jacobson} page 202, this implies $\Check{B}$
 $\in \langle g \rangle$.
 Hence $\Check{B}$ is of the form $a g + b e + c$ 
 for some $a, b, c \in \K$.
 Rewriting the Yang-Baxter equation for $\Check{B}$
 via the formulas above, 
 we obtain
 \begin{multline}
  ( \mu a^2 b + a b^2 ) X +
  ( - \mu a^2 c + a c^2 ) Y \\
  + \{ b^3 + \mu^2 a^2 b + (\lambda + 2 \mu) a b^2 + \lambda^{-1} \mu a^2 c
  + \zeta b^2 c + b c^2 + 2 \lambda^{-1} a b c \} Z
  = 0,  
 \end{multline}
 where
 \begin{equation}
  X = e_1 g_2 + g_2 e_1 - e_2 g_1 - g_1 e_2,
 \quad
  Y = g_1 - g_2,   
 \quad
  Z = e_1 - e_2.   
 \end{equation}
 Since $X, Y, Z$ are  linearly independent, we obtain
 three algebraic equations for $a$, $b$ and $c$.
 Solving these, we see that $\Check{B}$ is proportional to either
 $g$, $g^{-1} = g - \mu e + \mu$, $1 + \alpha e$ or $1$, where 
 $\alpha$ denotes a solution of 
 $x^2 + \zeta x + 1 = 0$.				
 Suppose $\Check{B} = \eta$ or $\eta (1 + \alpha e)$ 
 for some $\eta \in \K^{\times}$.
 Then using \eqref{mid(R)}, we obtain 
 \begin{equation}	
  \cal{B} (t_{12} t_{21},\, t_{22}) = \eta^2, 
 \quad
  \cal{B} (t_{21} t_{12},\, t_{22}) = 0.
 \end{equation}
 On the other hand, substituting $t_{21} \otimes t_{12}$
 into \eqref{RmXR}, we obtain 
 $t_{21} t_{12} = t_{12} t_{21}$, a contradiction.
 Therefore $\Check{B}$ is proportional to either
 $g$ or $g^{-1}$.
 This completes the proof of the theorem.
\end{pf} 			
The following lemma allows us to apply our general 
results developed in Sect. 6 to the Hopf closures.
\begin{lem}
 For each $q$ and $\eta$, $\K^N = \K \G^1$ is
 absolutely irreducible as
 a comodule of $\Hc ( \frak{A} ( \eta \Check{R}_q (X_l) ))$.
 In particular, $\Hc ( \frak{A} ( \eta \Check{R}_q (X_l) ))$
 is monogenerated.\\
 \end{lem} 
Since the proof of this lemma is quite similar to 
that of Lemma \ref{irrforqcg} below, we omit it.
Next, we determine the ribbon functionals of
the Hopf closure of 
$\frak{A} (\eta \Check{R}_q (X_l))$.  
We note that the following result immediately follows from 
Theorem \ref{Mcrit} and the formula \eqref{S2tij}
given below,
except for the case
$\sum_i  q^{2i - N - 1 - \sigma_i \nu} = 0$.
\begin{prop}
\label{clasribHc}
 For each $\eta \in \K^{\times}$, 
 $\Hc (\frak{A} ( \eta \Check{R}_q (X_l) ))$ 
 has exactly two \rom{(}resp. one\rom{)} modified ribbon functionals 
 $\cal{M}_{\pm}$ given by			
 \begin{equation}
 \label{M+-def}
 \cal{M}_{\pm} (t_{ij}) = 
 \pm \delta_{ij} 
 q^{2i - N - 1 - \sigma_i \nu}
 \end{equation}
 if $\mathrm{ch} \K \ne 2$
 \rom{(}resp. $\mathrm{ch} \K = 2$\rom{)}.
\end{prop}
\begin{pf}
 We will prove this result using Proposition \ref{ribestimate}
 and Theorem \ref{clasRibHcAw}.
 Using \eqref{S2tij} and \eqref{UXU-}, we obtain
 \begin{gather}
  (\pi (\cal{U}_1) \otimes \mathrm{id})
  \circ \rho \circ \pi (\cal{U}_1)^{-1} (u_j) =
  \sum_i u_i \otimes 
  q^{2(i -j) - (\sigma_i - \sigma_j) \nu}\,
  t_{ij} \\
  =
  (M \otimes \mathrm{id})
  \circ \rho \circ M^{-1} (u_j),
 \end{gather} 
 where  $M:= \mathrm{diag} 
 (q^{2i - N - 1 - \sigma_i \nu})_i$.
 This shows that $M \pi (\cal{U}_1)^{-1}$ commutes with
 the coaction of  $\Hc (\frak{A} ( \eta \Check{R}_q (X_l) ))$ 
 on $\K \G^1$. 
 Hence, it suffices to verify that
 \begin{equation}
 \label{piU1U2=M2} 
  \pi (\cal{U}_1 \cal{U}_2) = M^2.
 \end{equation}
 By Schur's lemma, we have  
 $\pi (\cal{U}_{\nu}) = c_{\nu} M$ for some constant 
 $c_{\nu} \in\K^{\times}$.
 Suppose $X = B, C$ or $D$. Using \eqref{SBCD}, 
 we compute 
 \begin{equation}
  \cal{U}_1 (t_{11}) = 
  \sum_{k=1}^{N} \epsilon_1 \epsilon_k
  q^{\overline{1} - \overline{k}} \eta \Check{R} 
  \left( k^{\prime} \frac[0pt]{1}{k} 1^{\prime} \right)
  = \eta \Check{R} \left( N \frac[0pt]{1}{1} N \right)
  = \eta q. 
 \end{equation}
 Using \eqref{S(U)} and \eqref{SBCD}, 
 we also obtain
 \begin{equation}
  \cal{U}_2^{-1} (t_{NN}) = \cal{U}_1 (t_{11}) 
  = \eta q.
 \end{equation}
 This proves $c_1 = \eta q^{N + \nu} = c_2^{-1}$,
 or \eqref{piU1U2=M2} for $X = B, C$, $D$.
 When $X = A$, \eqref{piU1U2=M2} is proved by computing 
 the Lyubashenko double of $\Check{R}_q (A_l)$
 explicitly. 
\end{pf}

Hereafter, we assume that $q^2 \ne -1$ and that 
$\lambda^{-1} \ne q^{-1}, -q$ if $X = B, C$ or $D$.
By \eqref{minpolyR}, this implies 
\begin{equation}
 \K^N \otimes \K^N = 
 \begin{cases}
  \mathrm{Ker} (\Check{R} - q^{-1}) \oplus
  \mathrm{Ker} (\Check{R} + q) 
    & (X = A)\\
  \mathrm{Ker} (\Check{R} - q^{-1}) \oplus
  \mathrm{Ker} (\Check{R} + q) \oplus
  \mathrm{Ker} (\Check{R} - \lambda^{-1})
  & (X = B, C, D)\\
 \end{cases}
\end{equation}
as $\frak{A} (\Check{R})$-comodules.
To give the definition of the quantized classical groups, 
we recall the definition of 
the (quantum)  determinant of 
$\frak{A} ( \Check{R})$.
Let $\Omega = \Omega (X_{l})$ be the following q-analogue of 
the exterior algebra:
\begin{equation}
 \Omega (X_{l}) =
 \begin{cases}
  T( \K^N ) / (\mathrm{Ker} (\Check{R} - q^{-1}))
  & (X = A, C) \\
  T( \K^N ) / (\mathrm{Im} (\Check{R} + q))
  & (X = B, D). 
 \end{cases}
\end{equation}
More explicitly, we have
\begin{equation}
 \Omega (A_{l}) = 
 \bigl\langle
 u_i \; (1 \leq i \leq N) \, 
 \bigm\vert \,
 u_i^2 = 0,\, q u_i u_j + u_j u_i = 0 \;
 (i < j)
 \bigr\rangle, 
\end{equation}
\begin{multline}
 \Omega (X_{l}) = 
 \bigl\langle
 u_i \; (1 \leq i \leq N) \,
 \bigm\vert \,
 u_i^2 = 0 \; (i \not= (N+1)/2), \\ 
 q u_i u_j + u_j u_i = 0 \;
 (i < j, i \not= j'), \; 
 \eta_i = 0 \;
 (1 \leq i \leq (N+1)/2)
 \bigr\rangle \\
 (X = B, C, D).
\end{multline}
Here for $X =$ $B$, $C$, $D$ and  
$1 \leq i \leq (N+1)/2$, we set
\begin{equation}
 \eta_i =
 \begin{cases}
 u_{i'} u_i +  u_i u_{i'} - 
 (q - q^{-1}) \sum_{j = 1}^{i-1} 
 q^{j - i + 1} u_j  u_{j'} &
 (X = B, D,\, i \leq l) \\
 u_{l+1} u_{l+1} -
 (q^{1/2} - q^{- 1/2}) \sum_{j = 1}^{l} 
 q^{j - l} u_j  u_{j'} &
 (X = B, \, i = l+1) \\
 u_{i'} u_i +  q^2 u_i u_{i'} +
 (q - q^{-1}) \sum_{j = i+1}^{l} 
 q^{j - i + 1} u_j  u_{j'} &
 (X = C,\, i \leq l).
 \end{cases}
\end{equation}
Then $\Omega (X_{l})$ becomes an 
$\frak{A} (\Check{R}_q (X_l))$-comodule algebra
via $u_j \mapsto \sum_i u_i \otimes t_{ij}$.
For $0 \leq r \leq N$, 
$\Omega_r := \sum_{i_1, \ldots, i_r} \K u_{i_1} \cdots u_{i_r}$
is a $\binom{N}{r}$-dimensional subcomodule of $\Omega$.
In particular, $\Omega_N$ $=$ $\K u_1 u_2 \ldots u_N$
is one-dimensional and determines the group-like element 
$\det \in \frak{A} (\Check{R})$ via the coaction
$u_1 \ldots u_N \mapsto 
u_1 \ldots u_N \otimes \det$.
For $X = B, C, D$, 
$\frak{A} (\Check{R})$ has another
group-like element $\mathrm{quad}$ which is determined by
its coaction on the one-dimensional comodule 
\begin{equation}
  \mathrm{Ker} (\Check{R} - \lambda^{-1}) =
 \K \sum_i \epsilon_i q^{\overline{i} + 1 / 2}
 u_i \otimes u_{i'}.
\end{equation}
By \cite{qcg} Proposition 5.4-5.5 and the universal 
mapping property of the Hopf closure and the localization
construction, we have
\begin{equation}
\label{HcA}
 \Hc ( \frak{A} ( \eta \Check{R}_q (A_l) )) \cong
 \frak{A} ( \eta \Check{R}_q (A_l) ) 
 [\mathrm{det}^{-1}]
 =:
 \mathrm{Fun} \left( GL_{q} (N) \right)_{\eta},
\end{equation}
\begin{align}
\label{HcBCD}
 \Hc ( \frak{A} ( \eta \Check{R}_q (X_l) )) 
 & \cong
 \frak{A} ( \eta \Check{R}_q (X_l) ) 
 [\mathrm{quad}^{-1}] \\
 & \cong 
 \frak{A} ( \eta \Check{R}_q (X_l) ) 
 [\mathrm{det}^{-1}]
 \quad (X = B, C, D).
\end{align}
The biideal $(\det - 1)$ becomes a CQT biideal of 
$\frak{A} (\eta \Check{R})$ if and only if
\begin{equation}
 \eta^N =
 \begin{cases}
  q & (X = A) \\
  1 & (X = B, C, D),
 \end{cases}
\end{equation} 
while $(\mathrm{quad} - 1)$
becomes a CQT biideal if and only if
$\eta = \pm 1$ (cf. \cite{gd}). 
%
%
%
%

Now we define the 
{\it function algebra of
the quantized classical groups}
(cf. \cite{TakeuchiMat}, \cite{qcg}, \cite{gd})
to be the CQT bialgebras given by
\begin{equation}
\label{SLdef}
 \mathrm{Fun} \left( SL_{q} (N) \right)_{\eta} =
 \frak{A} (\eta \Check{R}_q (A_l)) /
 (\mathrm{det} -1 )
 \quad (\eta^N = q),
\end{equation} 
\begin{multline}
\label{SOdef}
 \mathrm{Fun} \left( SO_{q} (N) \right)_{\eta} =
 \frak{A} (\eta \Check{R}_q (X_{l})) /
 (\mathrm{det} - 1, \mathrm{quad} - 1) \\
 (\eta = 1\, \mathrm{if}\, X = B\, \mathrm{and}\, 
 \eta = \pm 1\, \mathrm{if}\, X = D),
\end{multline} 
\begin{equation}
\label{Odef}
 \mathrm{Fun} \left( O_{q} (N) \right)_{\eta} =
 \frak{A} (\eta \Check{R}_q (X_{l})) /
 (\mathrm{quad} - 1) 
 \quad (\eta = \pm 1, X = B, D) 
\end{equation} 
\begin{equation}
\label{Spdef}
 \mathrm{Fun} \left( Sp_{q} (N) \right)_{\eta} =
 \frak{A} (\eta \Check{R}_q (C_{l})) /
 (\mathrm{quad} - 1)
 \quad (\eta = \pm 1). 
\end{equation} 
See \cite{qcg} for a justification of these definitions
in case $\K = \C$ and $q$ is transcendental over $\Bbb{Q}$.  
For $G_q = SL_q (N), SO_q (N)$, etc., 
we denote by $\mathrm{Fun} (G_q)$
the underlying bialgebra of $\mathrm{Fun} (G_q)_{\eta}$,
and by $\Rp_{\eta, G_q}$ the braiding of
$\mathrm{Fun} (G_q)_{\eta}$.  	
Each of these algebras has an antipode.
For example, the antipode of the algebras given in
\eqref{SOdef}-\eqref{Spdef}
is given by		
\begin{equation}
\label{SBCD}
 S(t_{ij}) = 
 \epsilon_i  \epsilon_j\,
 q^{\overline{i} - \overline{j}}\,
 t_{j'i'}.
\end{equation} 
The square of the antipode of the algebras given in
\eqref{HcA}, \eqref{HcBCD}, \eqref{SLdef}-\eqref{Spdef} 
is given by
\begin{equation}
\label{S2tij}
 S^2 (t_{ij}) = 
 q^{2(i -j) - (\sigma_i - \sigma_j) \nu}\,
 t_{ij}.
\end{equation} 
\begin{lem}
\label{irrforqcg}
 Let $F$ be either $\mathrm{Hc} (\frak{A} (\Check{R}))$ 
 or one of the algebras given in 
 \eqref{SLdef}-\eqref{Spdef}.
 Then each of the $F$-comodules $\K^N$, 
 $\mathrm{Ker} (\Check{R} - q^{-1})$ and  
 $\mathrm{Ker} (\Check{R} + q)$ 
 are absolutely irreducible.
 In particular, $F$ is monogenerated.\\
\end{lem} 	
\begin{pf}
 Since $\K$ is arbitrary, it suffices to show the irreducibility
 of these comodules. Here we give a proof for
 $W := \mathrm{Ker} (\Check{R}_q (C_l) + q)$.
 To simplify the computation, it is convenient to identify $W$
 with its image via the projection $(\K^N)^{\otimes 2} \to \Omega_2$.
 
 Following \cite{RTF}, 
 we define $K_i, E_i, F_i \in F^*$ $(1 \leq i \leq l)$ by 
 \begin{equation}
 \label{Kidef}
  K_i = \Rm_{\eta} (t_{ii},\, -), 
 \end{equation} 
 \begin{equation}
  E_i =
 \begin{cases}
  - \eta^{-1} (q - q^{-1})^{-1} 
  \Rp_{\eta} (-,\, t_{i+1\,\, i}) &
  (1 \leq i < l)\\
  - \eta^{-1} q^{-1} (q^2 - q^{-2})^{-1} 
  \Rp_{\eta} (-,\, t_{l+1\,\, l}) &
 (i = l), \\
 \end{cases}
 \end{equation} 
 \begin{equation}
  F_i =
 \begin{cases}
  \eta (q - q^{-1})^{-1} 
  \Rm_{\eta} (t_{i\,\, i+1},\,-) &
  (1 \leq i < l)\\
  \eta q (q^2 - q^{-2})^{-1} 
  \Rm_{\eta} (t_{l\,\, l+1},\,-) &
 (i = l). \\
 \end{cases}
 \end{equation} 
 Then these elements belong to the dual Hopf algebra
 $F^{\circ}$ (cf. \cite{LarsonTowber}) and satisfy 
 \begin{equation}
  \pi_{\K^N} (K_i) =
  \eta^{-1} \sum_{k=1}^{N} 
  q^{\delta_{ki} - \delta_{ki^{\prime}}} E_{kk}, 
 \end{equation} 
 \begin{equation}
  \pi_{\K^N} (E_i) = 
  E_{i\, i+1} - q E_{(i+1)^{\prime}\,  i^{\prime}}
  \quad (i < j),
 \quad
  \pi_{\K^N} (E_l) = 
  E_{l\,\, l+1}, 
 \end{equation} 
 \begin{equation}
  \pi_{\K^N} (F_i) = 
  E_{i+1\, i} - q^{-1} E_{i^{\prime}\, (i+1)^{\prime}}
  \quad (i < j),
 \quad
  \pi_{\K^N} (F_l) = 
  E_{l+1\,\, l}, 
 \end{equation} 
 \begin{equation}
  \Delta (K_i) = K_i \otimes K_i,
 \end{equation} 
 \begin{equation}
  \Delta (E_i) = E_i \otimes K_i^{-1} + K_{i + 1}^{-1} \otimes E_i,   
 \quad
  \Delta (F_i) = F_i \otimes K_{i + 1} + K_i  \otimes F_i,   
 \end{equation} 
 where $K_{l+1}$ is given by \eqref{Kidef}.

 As a $\langle K_i \rangle$-module, $W$
 is the direct sum
 of the mutually non-isomorphic, non-trivial comodules
 $\K u_i u_j$ $(j \ne i, i^{\prime})$ and the trivial comodule
 $T = \bigoplus_{i=1}^{l-1} \K (q u_{i} u_{i^{\prime}} -
 u_{i+1} u_{(i+1)^{\prime}})$. 
 Hence any non-zero subcomodule $M$ of $W$ contains a vector $v \ne 0$
 which belongs to one of these $\langle K_i \rangle$-modules.
 By verifying $T \cap (\bigcap_i \mathrm{Ker}\pi (E_i))$ $=$ $0$, 
 we see that $u_1 u_2 \in \K E_{i_1} \cdots E_{i_k} v$
 for some $i_1, \ldots, i_k$.
 Also, by verifying 
 $T = \sum_{i=1}^{l-1} \K F_i (u_i u_{(i+1)^{\prime}})$, 
 we see that $u_1 u_2$ generates $W$ as an
 $\langle F_i \rangle$-module.
 Thus, $W$ is irreducible as a
 $\langle K_i, E_i, F_i \rangle$-module, 
 and also, it is irreducible as an $F$-comodule.
\end{pf} 
\begin{thm}
\label{clasbrqcg}
\rom{(1)}		
 Any braidings of $\mathrm{Hc} (\frak{A} (\Check{R}_q (X_{l})))$
 are either of the form $\Rp_{\eta}$ 
 or of the form $(\Rm_{\eta})_{21}$, 
 where $\eta \in \K^{\times}$.\\
 \rom{(2)}		  
 Let $G_q$ be either $SL_q (N)$, $SO_q (N)$,  
 $O_q (N)$ or $Sp_q (N)$.
 Then, any braiding of $\mathrm{Fun} (G_q)$  	
 is either of the form $\Rp_{\eta, G_q}$ 
 or of the form $(\Rm_{\eta, G_q})_{21}$, where $\eta$ is 
 as in \eqref{SLdef}-\eqref{Spdef}.
\end{thm}
\begin{pf}
 Let $\Check{B}$ be a solution of the Yang-Baxter equation, 
 which corresponds to a braiding of 
 one of the above algebras $F$. 
 By Lemma \ref{irrforqcg}, 
 $\End_F ((\K^N)^{\otimes 2})$ is spanned by two or three projections 
 onto eigenspaces of $\Check{R}$, according to $X = A$ or $X = B, C, D$.
 By linear algebra, these projections are polynomials of $\Check{R}$.
 Therefore, we have  
 $\End_F ((\K^N)^{\otimes 2})$ $=$
 $\langle \Check{R} \rangle$.
 Hence, by the discussions in the proof of Theorem \ref{clasbrFRT}, 
 we see that $\Check{B}$ is proportional to either 
 $\Check{R}$ or $\Check{R}^{-1}$. 
 Thus this theorem follows from the result of \cite{gd} stated above.	
\end{pf} 
Let $F$ be one of the Hopf algebras given in \eqref{SLdef}-\eqref{Spdef}.
We define the cyclic group $\Gamma$ $=$ $\Gamma_F$
as follows:
\begin{equation}
 \Gamma =  
 \begin{cases}
  \Z / N \Z & ( G_q = SL_q(N))\\
  \Z / 2 \Z & ( G_q = O_q(N),\, SO_q (2 l),\, Sp_q (N))\\
  \{ 1 \} & ( G_q = SO_q (2 l + 1)).\\
 \end{cases}
\end{equation}
For $F$ $=$ $\frak{A} (\Check{R})$ and
$\mathrm{Hc} (\frak{A} (\Check{R}))$,  we also set
$\Gamma_F = \Z$.
Since $\det \in \frak{A}_N (\Check{R})$ and 
$\mathrm{quad} \in \frak{A}_2 (\Check{R})$,
the grading $\frak{A} (\Check{R})$ $=$ 
$\bigoplus_n \frak{A}_n (\Check{R})$ naturally induces
a $\Gamma$-grading of $F$ satisfying the properties 
stated in Proposition \ref{altbraiding} (2).
Now we can restate our classification theorems
for braidings as follows.		
\begin{cor}
\label{clasbrviaGamma}
 Let $F$ be one of the bialgebras treated in Theorem
 \ref{clasbrFRT} and Theorem \ref{clasbrqcg} and let
 $\Gamma$ be the cyclic group defined as above.
 Then, any braiding of $F$ is either of the form 
 $\Rp_{\chi}$ or of the form $(\Rm_{\chi})_{21}$, 
 where $\chi$ is as in Proposition \ref{altbraiding}.
 In particular, $\mathrm{MRib} (F)$ does not depend on
 the choice of the braiding of $F$ 
 \rom{(}cf. Proposition \ref{MRibaltbraiding}\rom{)}.
\end{cor}
Next, we give the classification theorem of the ribbon functionals 
for the algebras given in \eqref{SLdef}-\eqref{Spdef}.
\begin{lem}
 Let $\cal{M}_{\pm}$ be as in \eqref{M+-def}.
 Then we have
 \begin{equation}
  \cal{M}_{\pm} (\mathrm{det}) = (\pm 1)^N, 
 \quad 
  \cal{M}_{\pm} (\mathrm{quad}) = 1. 
 \end{equation}
\end{lem}
\begin{pf}
 We calculate
 \begin{multline}
  \cal{M}_{\pm} (\mathrm{det}) u_1 \cdots u_N =
  \pi_{\Omega_N}
  (\cal{M}_{\pm}) (u_1 \cdots u_N) =
  (\cal{M}_{\pm} u_1) \cdots 
  (\cal{M}_{\pm} u_N) \\
  = \prod_i \left( \pm q^{2i - N -1 - \sigma_i \nu} \right) 
  u_1 \cdots u_N 
  = (\pm 1)^N u_1 \cdots u_N.
  \qquad
 \end{multline} 
 The proof of the second formula is similar.
\end{pf} 
In view of the universal mapping property of the Hopf closure,
we see that we may replace $\frak{A} (\eta \Check{R})$ in
\eqref{SLdef}-\eqref{Spdef} with 
$\Hc( \frak{A} (\eta \Check{R}))$.
Hence, as an immediate consequence of Proposition \ref{clasquotient}
and the lemma above, 
we obtain the following.
\begin{thm}
\rom{(1)}
 Let $G_q$ be either  
 $SL_{q} (N)$, $SO_{q} (N)$ or $Sp_{q} (N)$,
 and let
 $p\!: \Hc (\frak{A} ( \eta \Check{R}_q (X_l) )) \to$ 
 $\mathrm{Fun} (G_q)_{\eta}$ denote the projection.
 Then we have		
 \begin{equation}
 \mathrm{MRib} \left(\mathrm{Fun} (G_q)_{\eta} \right) = 	
 \begin{cases}		
  \{ \cal{M}_+ \circ p,\, \cal{M}_- \circ p \} &
  (N \in 2 \Z) \\  
  \{ \cal{M}_+ \circ p \} &
  (N \in 1 + 2 \Z),
  \end{cases}		
 \end{equation}
 where $\eta$ is as in \eqref{SLdef}-\eqref{Spdef}. \\
\rom{(2)}
 We have
 \begin{equation}
 \mathrm{MRib} \left(
 \mathrm{Fun} \left( O_{q} (N) \right)_{\pm 1}
 \right) = 		
  \{ \cal{M}_+ \circ r,\, \cal{M}_- \circ r \}, 
 \end{equation}		
 where $r\!: \Hc( \frak{A} ( \pm \Check{R}_q (X_l) )) \to$ 
 $\mathrm{Fun} \left( O_{q} (N) \right)_{\pm 1}$
 denotes the projection.
\end{thm}

\section{SOS algebras}

%
%
%
Let $N \geq 2$ and $L \geq 2$ be integers. 
Let $\cal{C}$ be an $\C$-abelian semisimple rigid monoidal category  
whose simple objects $L_{\lambda}$ are indexed 
by the following set of partitions:
\begin{equation}
\label{Vdef}
 \V = {\V}_{NL}:= 	
 \bigl\{ \lambda = (\lambda_1, \ldots, \lambda_N) \in \Z^N
 \bigm|
 L \geq {\lambda}_1 \geq \dots \geq {\lambda}_{N} = 0
 \bigr\}. 
\end{equation} 
We say that $\cal{C}$ is an $SU(N)_L$-{\it category} if 
the structure constants of its Grothendiek ring
agree with the fusion rules $N^{\nu}_{\lambda \mu}$ 
of $SU(N)_L$-WZW models.
The $SU(N)_L$-categories play crucial roles to construct 
$SU(N)_L$-topological quantum field theories, or corresponding
invariants of 3-manifolds (cf. \cite{Turaev3mfd}).
It is known that two $SU(N)_L$-categories are equivalent to each 
other up to a ``twist'' of the associativity constraint
(cf. Kazhdan-Wenzl \cite{KazhdanWenzl}). 
	
In \cite{fut}, we have constructed a coribbon Hopf $\V_{NL}$-face algebra
$\frak{S} = \SS$ such that $\bold{Com}_{\frak{S}}^f$ is an $SU(N)_L$-category.  
In this section, we determine the braiding and the ribbon structure of 
$\bold{Com}_{\frak{S}}^f$ or equivalently, those of $\frak{S}$
(cf. Proposition \ref{correspalgcom}).

To begin with, we recall the definition of $SU(N)_L$-SOS model.
For each $1 \leq i \leq N$, we set
$\hat{i} = (\delta_{1i}, \ldots, \delta_{Ni}) \in \Z^N$.
For $m \geq 0$, we define the subset ${\G}^m$ of ${\V}^{m+1}$ by 
\begin{equation}
 {\G}^m = 
 {\V}^{m+1} \cap 
 \bigl\{ \p = (\lambda\, |\, i_1, \ldots, i_m)
 \bigm| \lambda\in \V,\, 
 1 \leq i_1, \ldots, i_m \leq N \bigr\},   
\end{equation}
where for $\lambda \in \Z^N$ and 
$1 \leq i_1, \ldots, i_m \leq N$,
we set 
\begin{equation}  
 (\lambda\, |\, i_1, \ldots, i_m) =
 (\lambda, \lambda + \hat{i}_1, \ldots, 
 \lambda + \hat{i}_1 + \cdots + \hat{i}_m),
\end{equation}
and we identify $(\lambda_1 +1, \cdots, \lambda_N + 1) \in \Z^N$ 
with $\lambda \in \V$.
Then $(\V,{\G}^1)$ defines an oriented graph $\G = {\G}_{N,L}$ 
and ${\G}^m$ is identified with the set of paths of $\G$ of length $m$.    
For $\p = (\lambda\, |\, i,j)$, 
we set 
${\p}^{\dag} = (\lambda\, |\, j,i)$ and
$d(\p) = {\lambda}_i - {\lambda}_j + j - i$. 
We define subsets
${\G}^2[\to]$, ${\G}^2[\;\downarrow\;]$ and
${\G}^2[\searrow]$ of ${\G}^2$ by				
\begin{gather}
 {\G}^2[\to]  =
 \bigl\{ \p \in {\G}^2 
 \bigm| \p^{\dag} = \p \bigr\},  
\quad 
 {\G}^2[\;\downarrow\;]  =
 \bigl\{ \p \in {\G}^2 
 \bigm| 
 \p^{\dag} \not\in \G^2 \bigr\},
 \\
 {\G}^2[\searrow] =  
 \bigl\{ \p \in {\G}^2 
 \bigm| \p \not= \p^{\dag} \in \G^2 
 \bigr\}.  
\end{gather} 
Let $t \in\C$ be a primitive $2(N+L)$-th root of $1$.
Let $\epsilon$ be either $1$ or $-1$
and $\zeta$ a nonzero parameter.
We define a star-triangular face model 
$(\G,w_{N,t,\epsilon}) = ({\G}_{N,L},w_{N,t, \epsilon, \zeta})$ 
by setting
\begin{equation}
 w_{N,t,\epsilon} \!
 \begin{bmatrix}
   \lambda           & \lambda + \hat{i} \\
   \lambda + \hat{i} & \lambda + \hat{i} + \hat{j}
 \end{bmatrix}  
 = - \zeta^{-1} t^{-d(\p)} \frac{1}{[d(\p)]},
\end{equation} 
\begin{equation}
 w_{N,t,\epsilon} \!
 \begin{bmatrix}
   \lambda           & \lambda + \hat{i} \\
   \lambda + \hat{j} & \lambda + \hat{i} + \hat{j}
 \end{bmatrix}  
 = \zeta^{-1} \epsilon \, \frac{[d(\p)-1]}{[d(\p)]}, 
\end{equation}
\begin{equation}
 w_{N,t,\epsilon} \!
 \begin{bmatrix}
   \lambda           & \lambda + \hat{k} \\
   \lambda + \hat{k} & \lambda + 2 \hat{k}
 \end{bmatrix}  
 = \zeta^{-1} t
\end{equation}
for each 
$\p = (\lambda\, |\, i,j) \in 
\G^2[\searrow] \amalg \G^2 [\,\downarrow\,]$
and 
$(\lambda\, |\, k,k) \in \G^2[\to]$, 
where $[n] = (t^n-t^{-n})/(t-t^{-1})$
for each $n \in \Z$.
We call $(\G,w_{N,t,\epsilon})$ {\it $SU(N)_L$-SOS model} 
(without spectral parameter) \cite{JMO}.
Now the $SU(N)_L$-SOS algebra $\SS$ is defined as the following
quotient of the FRT construction $\AS$:
\begin{equation}
 \SS:= \AS / (\det - 1),
\end{equation}
where the group-like element 
${\det} = \sum_{\lambda, \mu \in \V} \det \binom{\lambda}{\mu}$ 
of $\AS$ is defined by
\begin{equation}
 \label{detformula}
 \det \binom{\lambda}{\mu} =
 \frac{D(\mu)}{D(\lambda)}
 \sum_{\p \in \G^N_{\lambda \lambda}}
 (- \epsilon)^{\EuScript{L} (\p) + \EuScript{L} (\q)}
 e \binom{\p}{\q}.
\end{equation}
Here $\q$ denotes an arbitrary 
element of $\G^N_{\mu \mu}$, and 
$\EuScript{L} \!: \G^m \to \Z_{\geq 0}$
and $D( \lambda ) \in \C$ are given by
\begin{equation}
 \EuScript{L}
 (\lambda\, |\, i_1, \ldots, i_m)
 =
 \mathrm{Card} \{ (k,l) |
 1 \leq k < l \leq N, i_k < i_l \}, 
\end{equation}  
\begin{equation}
 D(\lambda) =
 \prod_{1 \leq i < j \leq N} 
 \frac{[d(\lambda\, |\, i,j)]}{[d(0\, |\, i,j)]}
 \quad (\lambda \in \V).
\end{equation}
The canonical braiding of $\frak{A} (w_{N,t,\epsilon, \zeta})$
induces the braiding $\Rp_{\zeta}$ of $\SS$ if and only if  
$\zeta$ satisfies
$\zeta^N = \epsilon^{N-1} t$. 
The face algebra $\SS$ has an antipode whose square is given by
\begin{equation}
 S^2 \left( e \binom{\p}{\q} \right) =
 \frac{D(\en (\p)) D(\st (\q))}{D(\st (\p)) D(\en (\q))}
 e \binom{\p}{\q}
 \quad (\p, \q \in \G^m, m \geq 0).
\end{equation}
The $\frak{S}$-comodule $\K \G^1$ is irreducible, while the 
$\frak{S}$-comodule $\K \G^2$ has the irreducible decomposition:
\begin{equation}
 \K \G^2 = \mathrm{Ker} (w_{N,t,\epsilon, \zeta} - \zeta^{-1} t) \oplus
 \mathrm{Ker} (w_{N,t,\epsilon, \zeta} + \zeta^{-1} t^{-1}).
\end{equation}
The proof of the following result is quite similar to
that of Theorem \ref{clasbrqcg} and
Takeuchi \cite{Takeuchicocycle} Lemma 2.4, 
hence we omit it.
\begin{thm}
 Any braiding of $\SS$ is either of the form $\Rp_{\zeta}$ or of the form
 $(\Rm_{\zeta})_{21}$, where $\zeta$ denotes a solution of
 $\zeta^N = \epsilon^{N-1} t$.
\end{thm}
Similarly to Corollary \ref{clasbrviaGamma}, we can rewrite
the result above in terms of the $\Z / N \Z$-grading of
$\SS$ induced by 
$\frak{A} (w_{N,t,\epsilon})$ 
$=$ $\bigoplus_n \frak{A}_n (w_{N,t,\epsilon})$. 
\begin{thm}
 When $N$ is odd, $\SS$ has exactly one ribbon functional.
 The corresponding modified ribbon functional
 $\cal{M}_+$ is given by
 \begin{equation}
 \label{M+SOS}
  \cal{M}_+ = \sumkl \frac{D(l)}{D(k)}\, \emak e_l.
 \end{equation}				
 When $N$ is even, $\SS$ has exactly two ribbon functionals.
 The corresponding modified ribbon functionals $\cal{M}_{\pm}$
 are given by \eqref{M+SOS} and 
 \begin{equation}
  \left\langle \cal{M}_-, 
  e\! \binom{\p}{\q} \right\rangle
  = \delta_{\p \q} (-1)^m
  \frac{D(\en (\p))}{D(\st (\p))}
  \quad (\p, \q \in \G^m, m \geq 0).
 \end{equation}				
 %
 %
 %
\end{thm}
\begin{pf}
 By Theorem \ref{Mcrit}, we have 
 $\mathrm{MRib} (\Hc ( \AS))$ $=$ $\{ \Phi (\pm M) \}$, 
 where $M \in GL (\K \G^1)$ is given by 
 $M \p = D(\en (\p)) D(\st (\p))^{-1} \p$ $(\p \in \G^1)$.
 Since
 \begin{align*}
 \left\langle \Phi (\pm M),\,  e \binom{\p}{\q} \right\rangle
  = & \left( \pm \delta_{\p_1 \q_1} 
  \frac{D(\en (\p_1))}{D(\st (\p_1))} \right) \cdots
  \left( \pm \delta_{\p_m \q_m} 
  \frac{D(\en (\p_m))}{D(\st (\p_m))} \right) \\
  = & (\pm 1)^m \delta_{\p \q} 
  \frac{D(\en (\p))}{D(\st (\p))}
 \end{align*}				
 for each $\p = (\p_1 \ldots \p_m)$, 
 $\q = (\q_1, \ldots \q_m) \in \G^m$,  
 we have 
 \begin{equation}		
  \langle \Phi (\pm M),\, \det - 1\rangle
  = \mathrm{Card} (\V) ( (\pm 1)^N - 1).
 \end{equation}				
 By Proposition \ref{clasquotient}, this proves the assertion.
\end{pf}
%

\end{document}